\documentclass[10pt,reqno]{amsart}
\usepackage[utf8]{inputenc} 
\usepackage[T1]{fontenc}    
\usepackage{lmodern}  

\usepackage{graphics,color}
\usepackage{amssymb, amsmath, amsthm, amscd}
\usepackage{latexsym, verbatim, graphicx, amsfonts}
\usepackage{hyperref}
\usepackage[mathscr]{euscript}
\usepackage{amsmath, amsthm, amssymb}
\usepackage{dsfont}
\usepackage{enumitem}
\usepackage{mathrsfs}
\usepackage{mathtools}
\usepackage{float}

\theoremstyle{plain}
\newtheorem{theorem}{Theorem}[section]

\newtheorem{corollary}[theorem]{Corollary}
\newtheorem{lemma}[theorem]{Lemma}
\newtheorem{proposition}[theorem]{Proposition}
\theoremstyle{definition}

\newtheorem{definition}[theorem]{Definition}
\newtheorem{example}[theorem]{Example}
\theoremstyle{remark}
\newtheorem{rem}[theorem]{Remark}

\newtheoremstyle{questionstyle}
  {6pt}   
  {6pt}   
  {\normalfont} 
  {}      
  {\bfseries} 
  {.}     
  {.5em}  
  {}      

\theoremstyle{questionstyle}
\newtheorem{question}[theorem]{Question}

\DeclareMathOperator*{\dist}{dist}
\DeclareMathOperator{\spn}{span}

\newcommand{\norm}[1]{\left\lVert#1\right\rVert}

\DeclareMathOperator{\incidence}{\mathbf{IA}}
\DeclareMathOperator{\opposite}{op}

\DeclareMathOperator{\identity}{Id}

\newcommand{\inset}[1]{\left[{#1}\right]^{\omega}}

\newcommand{\vertiii}[1]{{\left\vert\kern-0.25ex\left\vert\kern-0.25ex\left\vert #1 
    \right\vert\kern-0.25ex\right\vert\kern-0.25ex\right\vert}}

\font\sstext=ecss1000
\font\sssub=ecss1000 at 7pt
\font\sssubsub=ecss1000 at 5pt

\newfam\ssfam
\textfont\ssfam=\sstext
\scriptfont\ssfam=\sssub
\scriptscriptfont\ssfam=\sssubsub

\def\CH{\mathsf{CH}} 

\subjclass[2020]{47L10,
47L20,
46H10,
46E15 
(primary); 
46H40, 
46B26,
47B01,
54D80  
(secondary).}


\begin{document}

\title[$C(K)$-spaces with few operators relative to posets]{$C(K)$-spaces with few operators relative to posets}
\author[A. Acuaviva]{Antonio Acuaviva}
\address{School of Mathematical Sciences,
Fylde College,
Lancaster University,
LA1 4YF,
United Kingdom} \email{ahacua@gmail.com}

\date{\today}

\begin{abstract}
Extending a method developed by Koszmider and Laustsen for constructing $C(K)$-spaces we produce families of $C(K)$-spaces with few operators relative to a partially ordered set $\mathcal{P}$. Using these spaces, we construct new $C(K)$-spaces whose closed operator ideals can be completely classified. Additionally, we use these spaces to resolve some questions regarding automatic continuity.
\end{abstract}

\maketitle

\bigskip
\section{Introduction and main results}

The question of which kinds of operators must exist on a Banach space has a long and deep history. In the separable case, we find the seminal works of Gowers and Maurey \cite{gowers1993unconditional}. They constructed a space in which every operator is a strictly singular perturbation of a scalar multiple of the identity, providing a negative solution to the unconditional basic sequence problem. These ideas have had a profound impact on Banach space theory, yielding counterexamples to numerous open problems; see \cite{maurey2003banach} for a comprehensive survey of these developments. This cycle of ideas culminated in the landmark resolution of the ``scalar-plus-compact'' problem by Argyros and Haydon \cite{argyros2011hereditarily}, which has since inspired numerous new research directions. \\

However, the first construction of a Banach space with ``few operators'' was not achieved in the separable setting. Shelah \cite{shelah1978banach} constructed a non-separable Banach space on which every operator is a separable perturbation of a scalar multiple of the identity. Shelah's original construction relied on an additional set-theoretical axiom, $\Diamond$, but this was later removed in the work of Shelah and Stepr\=ans \cite{shelah1988banach}. This line of research was further developed by Argyros and Tolias \cite{argyros2004methods}, who constructed a dual space with this property, and by Wark \cite{wark2001non, wark2018non}, who produced a reflexive and a uniformly convex space with the aforementioned property. \\

All the previously mentioned Banach spaces have very complex definitions; in contrast, spaces of continuous functions over compact Hausdorff spaces $C(K)$ are among the simplest Banach spaces one can define. In this direction, in 2005, Koszmider  \cite{koszmider2005decompositions}, assuming the Continuum Hypothesis $(\CH)$, constructed a scattered locally compact Hausdorff space $K_\mathcal{A}$ such that the corresponding space $C_0(K_\mathcal{A})$ admits as few operators as possible, in the sense that every operator $T: C_0(K_\mathcal{A}) \to C_0(K_\mathcal{A})$ is of the form $T = \lambda \identity + S$, where $\lambda \in \mathbb{K}$ (the scalar field, either $\mathbb{R}$ or $\mathbb{C}$) and $S$ factors through $c_0$. The dependence on $(\CH)$ was later removed in the work of Koszmider and Laustsen \cite{koszmider2021banach}. Both arguments rely on the construction of a locally compact Hausdorff space induced by a carefully chosen almost disjoint family. \\

Let $\inset{\mathbb{N}}$ denote the set of all infinite subsets of the natural numbers. A family $\mathcal{A} \subseteq \inset{\mathbb{N}}$ is called \emph{almost disjoint} if for every distinct $A, B \in \mathcal{A}$, the intersection $A \cap B$ is finite. An almost disjoint family $\mathcal{A} \subseteq \inset{\mathbb{N}}$ induces a locally compact Hausdorff space $K_\mathcal{A}$, usually referred to as a \emph{$\Psi$-space}, \emph{Isbell--Mr\'owka space}, \emph{AU-compactum}, or \emph{Mr\'owka space}; see Definition \ref{def: almost-disjoint-family-top-space} for details.

These spaces were first introduced into Banach space theory by Johnson and Lindenstrauss, who used a renorming of $C_0(K_\mathcal{A})$-spaces to construct counterexamples concerning weakly compactly generated Banach spaces \cite{johnson1974some}.

Mr\'owka spaces have been used to construct several examples of exotic Banach spaces in recent years. Notably, Plebanek and Salguero-Alarc\'{o}n constructed an almost disjoint family $\mathcal{A}$ such that $C_0(K_\mathcal{A})$ provides a counterexample to the complemented subspace problem for $C(K)$-spaces \cite{plebanek2023complemented}. Remarkably, de Hevia, Martínez-Cervantes, Salguero-Alarc\'on, and Tradacete have subsequently shown that this space also serves as a counterexample to the complemented subspace problem for Banach lattices \cite{dehevia2025negativesolutioncomplementedsubspace}. We have recently extended these results by showing that a similar construction shows that the class of Banach lattices is not primary \cite{acuaviva2025class}. \\

As previously mentioned, these spaces were also used by Koszmider and Laustsen \cite{koszmider2021banach} to build a Banach space of continuous functions with as few operators as possible. It is within this context that our work is situated. We build on their ideas to construct new $C(K)$-spaces with interesting properties. In particular, these spaces allow us to resolve some of the questions left open in \cite{koszmider2021banach}. The following theorem summarises our main results in this direction, where $\wp(\mathfrak c)$ denotes the power set of the continuum, equipped with the order structure given by inclusion of sets.

\begin{theorem}\label{th: almost-disjoint-families-few-operators}
    Suppose that $\mathcal{P}$ is a partially ordered set which admits an order embedding into $\wp(\mathfrak c)$. Then there exist a collection $(\mathcal{A}_t)_{t \in \mathcal{P}}$ of uncountable, almost disjoint families $\mathcal{A}_t \subseteq \inset{\mathbb{N}}$ and a collection of isometric algebra homomorphisms $\iota_{s,t}: C_0(K_{\mathcal{A}_s}) \to C_0(K_{\mathcal{A}_t})$ for $s,t \in \mathcal{P}$ with $s \leq t$ such that:
    \begin{enumerate}
        \item $\iota_{s,s}$ is the identity on $C_0(K_{\mathcal{A}_s})$, and $\iota_{s,t} \iota_{r,s} = \iota_{r,t}$ for all $r \leq s \leq t$.
        \item Let $T: C_0(K_{\mathcal{A}_s}) \to C_0(K_{\mathcal{A}_t})$ be an operator for some $s,t \in \mathcal{P}$. Then
        \begin{itemize}[label=--]
            \item if $s \leq t$, then $T = \lambda \iota_{s,t} + S$, where $\lambda \in \mathbb{K}$ and $S$ is an operator that factors through $c_0$,
            \item if $s \not \leq t$, then $T$ factors through $c_0$.
        \end{itemize}
    \end{enumerate}
\end{theorem}

This theorem has many interesting consequences. In particular, every partially ordered set $\mathcal P$ with $|\mathcal P|\leq\mathfrak c$ order-embeds into $(\wp(\mathfrak c),\subseteq)$: the map $t\mapsto\{s\in\mathcal P:s\leq t\}$ embeds $\mathcal P$ into $(\wp(\mathcal P),\subseteq)$, and, after identifying $\mathcal P$ with a subset of $\mathfrak c$, this gives an embedding into $(\wp(\mathfrak c),\subseteq)$. 

Thus, for example, suppose that we let $\mathcal{P}$ be a set of cardinality continuum with the trivial order, that is, $s \leq t$ if and only if $s=t$. In that case, it gives continuum many non-isomorphic $C_0(K_\mathcal{A})$-spaces with few operators, thereby yielding a large collection of non-isomorphic $C(K)$-spaces for which the complemented subspaces and the closed operator ideals are completely classified, see \cite[Lemma 4]{koszmider2005decompositions} and \cite[Theorem 5.5]{kania2014ideal} respectively. Moreover, as an easy consequence, we recover the existence of $2^{\mathfrak c}$ distinct isomorphism classes of such spaces, a result previously obtained, by different methods, in \cite[Theorem 1.1(c)]{cabelloSanchez2020sailing}.

\begin{corollary}
There are $2^{\mathfrak c}$ pairwise non-isomorphic spaces of the form $C_0(K_{\mathcal A})$, where $\mathcal A$ is an almost disjoint family of cardinality $\mathfrak c$.
\end{corollary}

\begin{proof}
Choose an antichain $\mathcal P\subseteq \wp(\mathfrak c)$ of cardinality $2^{\mathfrak c}$. For example, identify $\mathfrak c$ with $\mathfrak c\times 2$, and for each function $\varepsilon:\mathfrak c\to 2$ define
\begin{equation*}
P_\varepsilon=\{(\xi,\varepsilon(\xi)):\xi<\mathfrak c\}.
\end{equation*}
Then $\mathcal P=\{P_\varepsilon:\varepsilon\in 2^{\mathfrak c}\}$ is an antichain in $(\wp(\mathfrak c),\subseteq)$ of cardinality $2^{\mathfrak c}$.

Applying Theorem~\ref{th: almost-disjoint-families-few-operators} to this partially ordered set, we obtain almost disjoint families $(\mathcal A_t)_{t\in\mathcal P}$ and spaces $C_0(K_{\mathcal A_t})$ for each $t\in\mathcal P$. We claim that these spaces are pairwise non-isomorphic.

Indeed, let $s,t\in\mathcal P$ with $s\neq t$. Since $\mathcal P$ is an antichain, we have $s\not\leq t$. Hence every operator from $C_0(K_{\mathcal A_s})$ into $C_0(K_{\mathcal A_t})$ factors through $c_0$, and in particular it has separable range. Thus, it cannot be an isomorphism, and therefore the spaces are not isomorphic.
\end{proof}

The true significance of the above theorem is that we can use the spaces $C_0(K_{\mathcal{A}})$ as building blocks to construct a plethora of $C(K)$-spaces displaying a rich variety of complemented subspaces and lattices of closed operator ideals. From this point onward, any reference to an operator ideal, or to lattices of operator ideals, will be understood to refer exclusively to closed operator ideals.

Analysing the closed ideals of $\mathscr{B}(C_0(K_\mathcal{A}))$ is intimately connected to identifying which unital Banach algebras can be isomorphic to the quotient algebra
\begin{equation*}
    \mathscr{B}(C_0(K_\mathcal{A}))/\mathscr{X}(C_0(K_\mathcal{A}))
\end{equation*}
for some almost disjoint family $\mathcal{A} \subseteq \inset{\mathbb{N}}$, where $\mathscr{X}$ denotes the ideal of operators with separable range. The following proposition, whose proof is essentially the same as that of \cite[Theorem 1.6]{acuaviva2025operators}, makes this connection explicit.

\begin{proposition}\label{prop: quotient-to-ideals}
    Let $\mathcal{A} \subseteq \inset{\mathbb{N}}$ be an almost disjoint family. The ideal $\mathscr{K}(C_0(K_\mathcal{A}))$ of compact operators is the smallest non-zero closed ideal of $\mathscr{B}(C_0(K_\mathcal{A}))$ and if $\varphi: \mathscr{B}(C_0(K_\mathcal{A})) \to \mathscr{B}(C_0(K_\mathcal{A}))/\mathscr{X}(C_0(K_\mathcal{A}))$ denotes the quotient map, then
    \begin{equation*}
        \mathscr{J} \mapsto \varphi^{-1}(\mathscr{J})
    \end{equation*}
    is a lattice isomorphism between the lattice of closed ideals of \break $\mathscr{B}(C_0(K_\mathcal{A})) / \mathscr{X}(C_0(K_\mathcal{A}))$ and the closed ideals of $\mathscr{B}(C_0(K_\mathcal{A}))$ strictly larger than $\mathscr{K}(C_0(K_\mathcal{A}))$.
\end{proposition}

\begin{rem}
    Let $\mathfrak{J}$ denote the lattice of closed ideals of \break $\mathscr{B}(C_0(K_\mathcal{A}))/\mathscr{X}(C_0(K_\mathcal{A}))$, then the previous proposition can be succinctly described by saying that the lattice of closed ideals of $\mathscr{B}(C_0(K_\mathcal{A}))$ is given by
    \begin{equation*}
        \{0\} \subsetneq \mathscr{K}(C_0(K_\mathcal{A})) \subsetneq \varphi^{-1}(\mathfrak{J}).
    \end{equation*}
\end{rem}

It follows that the problem of understanding closed ideals of \linebreak[4] $\mathscr{B}(C_0(K_\mathcal{A}))$ reduces to understanding the closed ideals of the quotient algebra $\mathscr{B}(C_0(K_\mathcal{A}))/\mathscr{X}(C_0(K_\mathcal{A}))$. Consequently, we focus on determining which unital Banach algebras can be realised as such quotients. The following terminology will allow us to state the main construction concisely.
\begin{itemize}
    \item Let $(\mathscr{C}_n)_{n \in \mathbb{N}}$ be a sequence of normed algebras. The \emph{normed direct product} of the sequence $(\mathscr{C}_n)_{n \in \mathbb{N}}$ is defined as
    \begin{equation*}
        \prod_{n \in \mathbb{N}} \mathscr{C}_n := \left\{a = (a_n)_{n \in \mathbb{N}}: a_n \in \mathscr{C}_n,  \hspace{5pt} \sup_{n \in \mathbb{N}} \norm{a_n} < \infty \right\},
    \end{equation*}
    equipped with the norm
    \begin{equation*}
        \norm{a} := \sup_{n \in \mathbb{N}} \norm{a_n},
    \end{equation*}
    and with addition and multiplication defined component-wise.
    \item Set
    \begin{align*}
        \mathfrak{A} = \{0\} 
        \cup \{\mathbb{M}_{n} : n \in \mathbb{N} \} 
        \cup \{\incidence(\mathcal{P}) : |\mathcal{P}| < \infty \} \cup \{\incidence(\mathcal{P}, c_0) : \mathcal{P} \text{ countable} \} 
        \cup \{\mathscr{B}(c_0)\},
    \end{align*}
    where $\mathbb{M}_n$ denotes the algebra of $n \times n$ matrices, equipped with the norm it inherits under the identification with $\mathscr{B}(\ell_\infty^n)$, $\incidence(\mathcal{P})$ is the incidence algebra induced by a finite poset $\mathcal{P}$, and $\incidence(\mathcal{P}, c_0)$ is the $c_0$-incidence algebra induced by a countable poset $\mathcal{P}$ (see Definitions \ref{def: inc-algebra} and \ref{def: E-inc-algebra} for details of this terminology).
\end{itemize}

\begin{theorem}\label{th: quotient-algebras}
Let $(\mathscr{C}_n)_{n \in \mathbb{N}} \subseteq \mathfrak{A}$. Then there exists an almost disjoint family $\mathcal{A} \subseteq \inset{\mathbb{N}}$ and a short exact sequence 
\begin{equation*}
    \{0\} \longrightarrow \mathscr{X}(C_0(K_\mathcal{A})) \longrightarrow \mathscr{B}(C_0(K_\mathcal{A})) \xlongrightarrow{\theta} \prod_{n \in \mathbb{N}} \mathscr{C}_n \longrightarrow \{0\}
\end{equation*}
which splits strongly. Furthermore, the continuous right inverse $\rho$ may be chosen to be an isometry, and $\norm{\theta}=1$ whenever $\prod_{n\in\mathbb N}\mathscr C_n\neq\{0\}$.
\end{theorem}

\begin{rem}
Here ``splits strongly'' is understood in the sense of Definition~\ref{def: strongly-split}; in particular, there is a continuous algebra homomorphism $\rho\colon \prod_{n \in \mathbb{N}}\mathscr{C}_n \to \mathscr{B}(C_0(K_\mathcal{A}))$ such that $\theta\circ\rho$ is the identity map on $\prod_{n \in \mathbb{N}}\mathscr{C}_n$. Equivalently, if zero factors are allowed, the one-factor case can be viewed as the product with $\mathscr C_1=\mathscr C$ and $\mathscr C_n=\{0\}$ for $n\geq2$.
\end{rem}

The previous theorem provides a partial answer to \cite[Question 46]{koszmider2021banach}. In particular, since we can obtain the finite product of matrix algebras, Wedderburn’s structure theorem \cite[Theorem 1.5.9]{dales2001banach} automatically gives the following.

\begin{corollary}
     Let $\mathscr{C}$ be a finite-dimensional, semi-simple complex algebra. Then there exists an almost disjoint family $\mathcal{A} \subseteq \inset{\mathbb{N}}$ such that $\mathscr{B}(C_0(K_\mathcal{A}))/\mathscr{X}(C_0(K_\mathcal{A}))$ is isomorphic to $\mathscr{C}$.
\end{corollary}

Using these families, we can prove the following result regarding automatic continuity.

\begin{theorem}\label{th:zfc-discontinuous-homomorphism}
There exists an almost disjoint family $\mathcal A\subseteq\inset{\mathbb N}$ such that there is a discontinuous surjective algebra homomorphism
\begin{equation*}
    \Psi\colon\mathscr B(C_0(K_\mathcal A))\longrightarrow\mathbb T_2(\mathbb K),
\end{equation*}
where $\mathbb T_2(\mathbb K)$ denotes the algebra of upper-triangular $2\times2$ scalar matrices. Consequently, there exist a compact scattered Hausdorff space $K$ of Cantor--Bendixson height at most three and a discontinuous surjective algebra homomorphism from $\mathscr B(C(K))$ onto $\mathbb T_2(\mathbb K)$.
\end{theorem}

The proof is deferred to Subsection~\ref{subsec:zfc-discontinuous}. This answers a question raised by the author \cite[Question 1.2]{acuaviva2024factorizations} and a question of Koszmider and Laustsen \cite[Question 47]{koszmider2021banach}. It also gives a genuinely new example and method, different from the one of Read \cite{read1989discontinuous}, of a Banach space $X$ such that $\mathscr{B}(X)$ admits a discontinuous homomorphism into a Banach algebra (without additional set-theoretic assumptions).

\begin{rem}
    Theorem~\ref{th:zfc-discontinuous-homomorphism} does not control derivations from the same operator algebra. Following the arguments in \cite{dales1994homomorphisms}, we believe that, under $(\CH)$, it should be possible to construct an almost disjoint family $\mathcal{A}$ such that every derivation from $\mathscr{B}(C_0(K_\mathcal{A}))$ is continuous while the algebra admits a discontinuous homomorphism into a Banach algebra. Establishing these two properties simultaneously would require a substantially longer argument.
\end{rem}

We want to highlight that there are extremely few $C(K)$-spaces for which the lattice of closed operator ideals has been fully classified. In fact, even the number of such ideals is known only for the spaces listed below. To the best of our knowledge, these are:
\begin{itemize}
    \item $c_0(\Gamma)$ for any set $\Gamma$; in this case, the lattice of closed operator ideals was classified by Gohberg, Markus and Feldman \cite{gohberg1967normally} in the case $\Gamma = \mathbb{N}$ and extended to uncountable sets by Daws \cite{daws2006closed}. This lattice is well-ordered and consists of the ideals of $\kappa$-compact operators. An alternative description of this ideal lattice was also given by Johnson, Kania and Schechtman \cite[Theorem 1.5]{johnson-kania-schechtman2016}.
    \item $C_0(K_\mathcal{A})$ where $K_\mathcal{A}$ is the space constructed by Koszmider and Laustsen \cite{koszmider2021banach} so that $C_0(K_\mathcal{A})$ admits few operators. In this case, Kania and Kochanek \cite[Theorem 5.5]{kania2014ideal}, showed that the lattice of closed ideals is
    \begin{equation*}
        \{0\} \subsetneq \mathscr{K}(C_0(K_\mathcal{A})) \subsetneq \mathscr{X}(C_0(K_\mathcal{A})) \subsetneq \mathscr{B}(C_0(K_\mathcal{A})).
    \end{equation*}
    A similar conclusion holds for $C_0(K_\mathcal{A})^n \cong C_0(K_\mathcal{A} \times \{1, \dots, n\})$ for any $n \in \mathbb{N}$.
    \item $C_0(\omega \times K_\mathcal{A})$ where again $K_\mathcal{A}$ is the space built by Koszmider and Laustsen. We recently showed that the lattice of closed operator ideals consists of five linearly ordered ideals \cite[Theorem 1.6]{acuaviva2025operators}.
\end{itemize}

It is interesting to note that all the previously known examples have a linearly ordered lattice and thus a unique maximal ideal. Applying Theorem \ref{th: quotient-algebras}, together with Proposition \ref{prop: quotient-to-ideals}, allows us to produce multiple $C(K)$-spaces whose lattices of closed operator ideals behave very differently, using the fact that the ideal structure of incidence algebras is essentially known \cite[Section 3]{doubilet1972foundations}.

In particular, we can obtain ideal lattices that are not linearly ordered and with many maximal ideals. In fact, for each $n \in \mathbb{N}$, the $n \times n$ upper triangular matrices have precisely $n$ maximal ideals, and thus we can find an almost disjoint family $\mathcal{A}$ such that $\mathscr{B}(C_0(K_{\mathcal{A}}))$ has $n$ maximal ideals. We refer to Subsection \ref{subsec:possible-quotient-algebras} for a more detailed account of possible ideals.

Lastly, we would like to emphasise that, to obtain spaces of operators on $C(K)$-spaces with more than one maximal ideal, one must consider non-separable $C(K)$-spaces. We refer to Section \ref{sec:separable-case} for details, where we show the following. 

\begin{proposition}\label{prop:separable}
    Let $C(K)$ be separable. Then $\mathscr{B}(C(K))$ has a unique maximal ideal. Moreover, there exists some $n \in \{2,3\}$ such that the space of operators on any separable $C(K)$-space never has exactly $n$ non-trivial closed ideals.
\end{proposition}

In contrast, in the non-separable setting, as previously mentioned, a result of Daws \cite{daws2006closed} shows that, for any $n \in \mathbb{N}$, the algebra $\mathscr{B}(c_0(\aleph_{n-1}))$ admits exactly $n$ non-trivial closed ideals. This motivates us to pose the following question.

\begin{question}\label{quest:numberideals}
For which $n \in \mathbb{N}$, $n \geq 1$, does there exist an almost disjoint family $\mathcal{A} \subseteq \inset{\mathbb{N}}$ such that the operator algebra $\mathscr{B}(C_0(K_\mathcal{A}))$ contains exactly $n$ non-trivial closed ideals?
\end{question}

We note that, by Theorem~\ref{th: quotient-algebras}, a significant subset of the natural numbers can be realised as the number of non-trivial closed ideals of algebras of the form $\mathscr{B}(C_0(K_\mathcal{A}))$. This is explored further in Subsection~\ref{subsec:possible-quotient-algebras}, where we give a general procedure for producing such numbers from finite incidence algebras and finite products of admissible quotient algebras.

We conjecture that, in fact, every natural number arises in this way, but a new procedure to build specific algebras as quotients may be needed.

\bigskip
\section{Organisation and notation}

All normed spaces and algebras are over the scalar field $\mathbb{K}$, either the real or complex numbers, and we adhere to standard notational conventions. For a compact Hausdorff space $\Omega$, we denote by $C(\Omega)$ the space of continuous functions $f: \Omega \to \mathbb{K}$. Similarly, for a locally compact Hausdorff space $K$, we denote by $C_0(K)$ the space of continuous functions on $K$ that \emph{vanish at infinity}, that is, continuous functions $f: K \to \mathbb{K}$ such that the set $\{k \in K: |f(k)| \geq \varepsilon \}$ is compact for each $\varepsilon > 0$.

The term \emph{operator} will refer to a bounded linear map between normed spaces. For two Banach spaces $X$ and $Y$, we denote by $\mathscr{B}(X, Y)$ the space of operators between $X$ and $Y$. In the case $X = Y$, we simply write $\mathscr{B}(X)$. When referring to ideals, we will always mean \emph{two-sided closed} ideals.

Furthermore, whenever we consider the finite direct sum of $n$ Banach spaces, it is understood implicitly that the sum is equipped with the $\ell_\infty^n$ norm. We will denote by $\inset{\mathbb N}$ the set of all infinite subsets of $\mathbb N$. Finally, given a sequence $(a_n)_{n \in \mathbb{N}}$ and a set $A \in \inset{\mathbb{N}}$, we write $\lim_{n \in A} a_n$ to mean $\lim_{k \to \infty} a_{n_k}$, where $(n_k)_{k \in \mathbb{N}}$ is the increasing enumeration of $A$.

We recall the following standard terminology for Banach algebras.

\begin{definition}
A sequence of Banach algebra homomorphisms
\begin{equation*}
    \{0\}\longrightarrow \mathscr I \xlongrightarrow{\jmath} \mathscr B \xlongrightarrow{\theta} \mathscr C \longrightarrow \{0\}
\end{equation*}
is called \emph{short exact} if $\jmath$ is injective, $\theta$ is surjective and $\ker \theta=\jmath(\mathscr I)$. In the applications below, $\jmath$ will be the inclusion map, and we shall identify $\mathscr I$ with its image in $\mathscr B$.
\end{definition}

\begin{definition}\label{def: strongly-split}
A short exact sequence
\begin{equation*}
    \{0\}\longrightarrow \mathscr I \longrightarrow \mathscr B \xlongrightarrow{\theta} \mathscr C \longrightarrow \{0\}
\end{equation*}
is called \emph{strongly split} if there is a continuous algebra homomorphism $\rho\colon \mathscr C\to \mathscr B$ such that $\theta\circ \rho$ is the identity map on $\mathscr C$. In this case, $\rho$ is called a \emph{continuous right inverse} for $\theta$.
\end{definition}

More specialised notation will be introduced as and when needed.

We give a brief overview of the structure of this paper. In Section~\ref{sec:separable-case}, we discuss the separable $C(K)$ case and show that separable $C(K)$-spaces always have a unique maximal ideal of operators. This explains why the later examples must be non-separable if one wants several maximal ideals. In Section~\ref{sec: preliminaries}, we introduce the notation and preliminary results needed for the proof of Theorem~\ref{th: almost-disjoint-families-few-operators}. In Subsection~\ref{sub-sec: topological-combination} we begin by recalling the definition of the topological spaces $K_\mathcal{A}$ induced by an almost disjoint family $\mathcal{A} \subseteq \inset{\mathbb{N}}$, as well as some basic properties concerning separability. We then explore how to combine different $C_0(K_\mathcal{A})$-spaces to construct new spaces of this form, and explain how the inclusion $\mathcal{A} \subseteq \mathcal{B} \subseteq \inset{\mathbb{N}}$ of almost disjoint families affects the relationship between the spaces $C_0(K_\mathcal{A})$ and $C_0(K_\mathcal{B})$.

In Subsection~\ref{subsec: matrices}, we describe how infinite matrices can be interpreted as acting between the spaces $C_0(K_\mathcal{A})$ and $C_0(K_\mathcal{B})$. To do this, we adapt several definitions and results originally developed by Koszmider and Laustsen \cite{koszmider2021banach} to the setting in which the domain and codomain differ. Readers already familiar with their work may wish to skip or only briefly skim this subsection, referring back to it as and when needed.

Lastly, Subsection~\ref{sub-sec: operators-and-matrices} explains how to represent an operator $T: C_0(K_{\mathcal{A}}) \to C_0(K_{\mathcal{B}})$ as an $\mathbb{N} \times \mathbb{N}$ matrix, and how this representation relates to the matrix theory discussed earlier. Readers familiar with the work of Koszmider and Laustsen will recognise most of the results presented here as natural extensions of ideas already developed in \cite{koszmider2021banach}. Section~\ref{sec: proof-main-theorem} presents the proof of Theorem~\ref{th: almost-disjoint-families-few-operators}.

Finally, Section~\ref{sec: operator-ideals} introduces the concepts of incidence algebras and $c_0$-incidence algebras, and explains how these algebras can be realised as quotients by the ideal of separable-range operators. Section~\ref{sec: proof-theorem-2} contains the proof of Theorem~\ref{th: quotient-algebras}. The paper concludes with Section~\ref{sec:applications}. Subsection~\ref{subsec:zfc-discontinuous} proves the ZFC automatic-continuity result stated in Theorem~\ref{th:zfc-discontinuous-homomorphism}, while Subsection~\ref{subsec:possible-quotient-algebras} discusses possible quotient algebras and ideal counts and records the multiplicative nature of the construction.

\bigskip
\section{Operator ideals for separable \texorpdfstring{$C(K)$}{C(K)}-spaces}\label{sec:separable-case}

This section explains why the examples constructed later must be non-separable if one wants several maximal ideals. It is easy to build non-separable $C(K)$-spaces whose operator ideal lattice is not linearly ordered. For instance, since the algebras of operators on $C([0,1])$ and on $C(\omega_1)$ each have a unique maximal ideal, one can verify that the algebra of operators on $C([0,1])\oplus C(\omega_1)$ has precisely two maximal ideals. However, a complete classification of the closed operator ideals on spaces of this kind appears to be out of reach. By contrast, in the separable $C(K)$ setting there is always a unique maximal ideal, as we will see.

We shall use the following notation. For Banach spaces $X$ and $Y$, define
\begin{equation*}
    \mathcal{M}_Y(X)=\{T\in\mathscr{B}(X): I_Y \text{ does not factor through } T\}.
\end{equation*}
When $X=Y$, this agrees with the notation introduced by Dosev and Johnson \cite{dosev2010commutators}, namely $\mathcal{M}_X(X)=\mathcal{M}_X$.
We can now make the following observation, which appears to have gone unnoticed in the literature.

\begin{proposition}\label{prop: maximal-ideal-alpha}
For every infinite ordinal $\alpha$, the algebra of operators on $C(\alpha)$ has a unique maximal ideal. In particular, if $C(\alpha) \not \sim C(\xi \cdot n)$ for any uncountable regular cardinal $\xi$ and $n \in \mathbb{N}, n \geq 2$ then the set $\mathcal{M}_{C(\alpha)}(C(\alpha))$ is the unique maximal ideal of $\mathscr{B}(C(\alpha))$.
\end{proposition}
\begin{proof}
    First, assume that $C(\alpha)$ is not isomorphic to $C(\xi \cdot n)$ for any uncountable regular cardinal $\xi$ and $n \in \mathbb{N}$, $n \geq 2$. In this case,  an inspection of the proof of the primariness of $C(\alpha)$ by Benyamini and Alspach \cite{alspach1977primariness} shows that if $T, S \in \mathscr{B}(C(\alpha))$ are such that $T + S = I_{C(\alpha)}$, then one of $T$ or $S$ factors the identity on $C(\alpha)$ (in the case that $\alpha$ is countable, this can also be seen from Wolfe's proof of primariness \cite{wolfe1982c}). 
    
    Indeed, for Cases~I and~II in the proof of \cite[Theorem 1]{alspach1977primariness}, we can still invoke \cite[Proposition 2]{alspach1977primariness} for any operator $T$, not necessarily a projection. The assumption $T + S = I_{C(\alpha)}$ is precisely what we need for arguing that we can assume $X$ is $(c, 1/10)$-preserved with $c \geq 1/2$ by either $T$ or $S$. After this, we can invoke \cite[Lemma 1.2(b) and (c)]{alspach1977primariness} to conclude that the operator, either $T$ or $S$, fixes a complemented copy of $C_0(\alpha)$ and thus it factors the identity $I_{C(\alpha)}$.

    Cases~III and~IV follow in the same way, noting that the fact that $P_Y$ and $P_Z$ are projections is never actually used; the only essential point is that $P_Y + P_Z = I_{C(\alpha)}$, which ensures that the operators $(c, 1/20)$-preserve the spaces $X_\delta$ for some $c \geq 1/2$.
    
    Equivalently, the preceding argument shows that for every operator $T \in \mathscr{B}(C(\alpha))$, we have
    \begin{equation*}
        T \notin \mathcal{M}_{C(\alpha)}(C(\alpha)) 
        \quad \text{or} \quad 
        I_{C(\alpha)} - T \notin \mathcal{M}_{C(\alpha)}(C(\alpha)).
    \end{equation*}
    It then follows from \cite[Proposition~5.1]{dosev2010commutators} that 
    $\mathcal{M}_{C(\alpha)}(C(\alpha))$ is the unique maximal ideal in 
    $\mathscr{B}(C(\alpha))$. 
    
    Suppose now that $C(\alpha) \sim C(\xi \cdot n) \cong \bigoplus_{j=1}^n C(\xi)$ for some uncountable regular cardinal $\xi$ and integer $n \ge 2$. By the previous argument, $\mathscr{B}(C(\xi))$ admits a unique maximal ideal. Elementary arguments using the matrix representation of operators on direct sums then yield the uniqueness of a maximal ideal in this case.
    
    By the classification theorem for spaces of continuous functions on ordinals \break —originally proved independently by Gul'ko \cite{gul1975isomorphic} and Kislyakov \cite{kislyakov1975classification}, and stated in the specific form used here in \cite{alspach1977primariness}—these are all the possible cases up to isomorphism.
\end{proof}

\begin{rem}
    The case $C(\alpha) \sim C(\omega_1)$, or more generally $C(\alpha) \sim C(\omega_\eta)$ where $\eta$ is any ordinal such that $\omega_\eta$ is a regular cardinal, was already shown by Kania and Laustsen \cite{kania2012uniqueness}, while the case $C(\alpha) \sim C(\beta)$ for $\beta$ a countable epsilon number, that is, a countable ordinal satisfying $\omega^\beta=\beta$, was shown by Philip A. H. Brooker (unpublished).
\end{rem}

We now pass from ordinal compacta to arbitrary separable $C(K)$-spaces.

\begin{corollary}\label{cor:uniqueseparable}
    Let $K$ be a compact metric space. Then $\mathscr{B}(C(K))$ admits a unique maximal ideal.
\end{corollary}
\begin{proof}
    If $K$ is finite, the result is trivial. Whenever $K$ is countable and infinite, the classification theorem of Bessaga and Pe{\l}czy{\'n}ski \cite{bessaga1960spaces} gives $C(K) \sim C(\alpha)$ for some infinite countable ordinal $\alpha$. The result now follows from Proposition \ref{prop: maximal-ideal-alpha}.

    If $K$ is uncountable, Milutin's theorem \cite{milutin1966isomorphisms} gives $C(K) \sim C([0,1])$. That $\mathscr{B}(C[0,1])$ has a unique maximal ideal follows by combining results of Pe{\l}czy{\'n}ski \cite[Theorem 1]{pelczynski1968c} and Rosenthal \cite[Theorem 1]{rosenthal1972factors} (see, for example, \cite[Example 3.5]{brooker2012asplund}).
\end{proof}

Given a natural number $n \in \mathbb{N}$, one may ask whether there exists a $C(K)$-space for which $\mathscr{B}(C(K))$ has exactly $n$ non-trivial closed ideals. We will show that this cannot be achieved in the separable setting. 

Observe that whenever $\mathcal{M}_Y(Y)$ is an ideal (that is, when it is closed under addition), then so is $\mathcal{M}_Y(X)$ for any Banach space $X$. In particular, by Proposition \ref{prop: maximal-ideal-alpha}, for any Banach space $X$ and any infinite countable ordinal $\alpha$, $\mathcal{M}_{C(\alpha)}(X)$ is an operator ideal.

\begin{proposition}\label{prop:separable1}
The algebra of operators on $C(\omega^{\omega^2})$ has at least four non-trivial closed ideals. Consequently, there exists some $n \in \{2,3\}$ such that the space of operators on any separable $C(K)$-space never has exactly $n$ non-trivial closed ideals.
\end{proposition}
\begin{proof}
    For the first part, we exhibit four different closed ideals in $\mathscr{B}(C(\omega^{\omega^2}))$. Recall that, for each ordinal $\alpha$, $\mathscr{SZ}_{\alpha}$, the collection of operators having Szlenk index at most $\omega^\alpha$, forms an operator ideal in the sense of Pietsch \cite[Theorem 2.2]{brooker2012asplund}.

We claim that the inclusions
\begin{equation*}
    \mathscr{K}(C(\omega^{\omega^2})) \subsetneq \mathscr{SZ}_{1}(C(\omega^{\omega^2})) \subsetneq \mathscr{SZ}_{2}(C(\omega^{\omega^2})) \subsetneq \mathcal{M}_{C(\omega^{\omega^2})}(C(\omega^{\omega^2}))
    \end{equation*}
    are all strict, and thus define four distinct non-trivial ideals in $\mathscr{B}(C(\omega^{\omega^2}))$. The first two strict inclusions follow from the standard computation of Szlenk indices for countable ordinal $C(K)$-spaces, together with the fact that $C(\omega^{\omega^2})$ contains complemented copies of $c_0$ and of $C(\omega^\omega)$. Thus one obtains a non-compact operator of Szlenk index at most $\omega$, and an operator of Szlenk index $\omega^2$.
    
    It remains to show that
    \begin{equation*}
        \mathcal{M}_{C(\omega^{\omega^2})}(C(\omega^{\omega^2})) \nsubseteq \mathscr{SZ}_{2}(C(\omega^{\omega^2})).
    \end{equation*}
    For this, observe that Alspach \cite{alspach1981c} constructed a surjective operator $T\colon C(\omega^{\omega^2}) \to C(\omega^{\omega^2})$ which does not preserve any copy of $C(\omega^{\omega^2})$. In particular, $T \in \mathcal{M}_{C(\omega^{\omega^2})}(C(\omega^{\omega^2}))$. On the other hand, since $T$ is onto and the Szlenk index of $C(\omega^{\omega^2})$ is $\omega^3$, we have $T \notin \mathscr{SZ}_{2}(C(\omega^{\omega^2}))$. This proves the first part. The second part is immediate.
\end{proof}

Proposition \ref{prop:separable}, as stated in the introduction, now follows by combining Corollary \ref{cor:uniqueseparable} and Proposition \ref{prop:separable1}.  
\bigskip
\section{Preliminary results}\label{sec: preliminaries}

\subsection{Structure and combination of \texorpdfstring{$C_0(K_\mathcal{A})$}{C0(KA)}-spaces.}\label{sub-sec: topological-combination} We briefly review the structure of $C_0(K_\mathcal{A})$-spaces and examine possible ways to combine them. We begin by recalling how an almost disjoint family $\mathcal{A} \subseteq \inset{\mathbb{N}}$ induces a scattered locally compact Hausdorff space $K_{\mathcal{A}}$.

\begin{definition}\label{def: almost-disjoint-family-top-space}
    Let $\mathcal{A} \subseteq \inset{\mathbb{N}}$ be an almost disjoint family. We define the topological space $K_\mathcal{A}$ consisting of distinct points $\{x_n: n \in \mathbb{N}\} \cup \{y_A: A \in \mathcal{A}  \}$, where $x_n$ is isolated for each $n \in \mathbb{N}$ and the sets
    \begin{equation*}
        U(A, F) = \{x_n: n \in A \setminus F \} \cup \{y_A\}
    \end{equation*}
    for finite $F \subseteq \mathbb{N}$ form a neighbourhood basis at each point $y_A$ for $A \in \mathcal{A}$. We will write $U(A)$ for $U(A, \emptyset)$.
\end{definition}

We refer the reader to \cite{hernandez2018topology} and \cite{hrusak2014almost} for a deeper discussion about these topological spaces. We note that the space $K_\mathcal{A}$ has finite Can\-tor–\-Ben\-dix\-son rank. Specifically, if $K'$ denotes the Cantor–Bendixson derivative of a topological space $K$, we have
\begin{equation*}
    K_\mathcal{A}' = \{y_A: A \in \mathcal{A}  \}, \text{ }K''_\mathcal{A} = \emptyset \text{ and } K_\mathcal{A} \setminus K_\mathcal{A}' = \{x_n: n \in \mathbb{N}\}.
\end{equation*}

We introduce a few concepts related to separability in $C_0(K_\mathcal{A})$.

\begin{definition}[{\cite[Definition 2]{koszmider2021banach}}]
    Let $\mathcal{A} \subseteq \inset{\mathbb{N}}$ be an almost disjoint family. For $f \in C_0(K_\mathcal{A})$ and $\mathcal{X} \subseteq C_0(K_\mathcal{A})$, we define
    \begin{equation*}
        s(f) = \{A \in \mathcal{A}: f(y_A) \not = 0\}, \hspace{5pt} s(\mathcal{X}) = \bigcup \{s(f): f \in \mathcal{X}\}.
    \end{equation*}
\end{definition}

\begin{lemma}[{\cite[Lemma 21]{koszmider2021banach}}]\label{lmm: s-countable}
     Suppose that $\mathcal{A} \subseteq \inset{\mathbb{N}}$ is an almost disjoint family. A closed subspace $\mathcal{X}$ of $C_0(K_\mathcal{A})$ is separable if and only if $s(\mathcal{X})$ is countable.
\end{lemma}

We will also need the following lemma, whose proof can be found in \cite[Lemma 3]{koszmider2005decompositions} or \cite[Example 2c]{johnson1974some}.

\begin{lemma}[{\cite[Lemma 22]{koszmider2021banach}}]\label{lmm: countable-to-fact-c0}
     Suppose that $\mathcal{A} \subseteq \inset{\mathbb{N}}$ is an almost disjoint family and $\mathcal{X} \subseteq C_0(K_\mathcal{A})$ is a separable subspace. Then there exists a closed subspace $\mathcal{Y} \subseteq C_0(K_\mathcal{A})$ such that $\mathcal{X} \subseteq \mathcal{Y}$ and $\mathcal{Y} \sim c_0$.
\end{lemma}

Our next result shows how spaces induced by almost disjoint families combine. In particular, it explains how a $c_0$-sum of $C_0(K_\mathcal{A})$-spaces can be seen as another $C_0(K_\mathcal{A})$-space.

\begin{lemma}\label{lmm: combine-pairwise-disjoint-families}
    Suppose that $(\mathcal{A}_n)_{n \in \mathbb{N}}$ is a sequence of almost disjoint families $\mathcal{A}_n \subseteq \inset{\mathbb{N}}$. Then there exists an almost disjoint family $\mathcal{B}$ such that $K_\mathcal{B}$ is homeomorphic to the topological sum $\coprod_{n \in \mathbb{N}} K_{\mathcal{A}_n}$.
    In particular, we have
  \begin{equation*}
        C_0(K_\mathcal{B}) \cong C_0 \left(\coprod_{n \in \mathbb{N}} K_{\mathcal{A}_n} \right) \cong \left(\bigoplus_{n \in \mathbb{N} } C_0(K_{\mathcal{A}_n})\right)_{c_0}.
    \end{equation*}
\end{lemma}
\begin{proof}
    Partition $\mathbb{N}$ into countably many disjoint infinite subsets $\mathbb{N} = \bigsqcup_{n \in \mathbb{N}} C_n$, $C_n \in \inset{\mathbb{N}}$, and enumerate $C_n = \{c_{n,m}: m \in \mathbb{N}\}$ for each $n \in \mathbb{N}$.

    For each $n \in \mathbb{N}$ and each $A \in \mathcal{A}_n$ define $B(n,A) = \{c_{n,m}: m \in A\}$. It is straightforward to check that the family $\mathcal{B} = \{B(n, A): n \in \mathbb{N}, A \in \mathcal{A}_n\}$ works.
\end{proof}

We finally explain how the inclusion of almost disjoint families affects the relation between the induced $C_0(K_\mathcal{A})$-spaces. In particular, if $\mathcal{A}$ is a subset of an almost disjoint family $\mathcal{B}$, it will allow us to isometrically identify $C_0(K_\mathcal{A})$ as a subalgebra of $C_0(K_\mathcal{B})$ in a canonical way. 

This turns out to be a consequence of a more general ``folklore'' result, communicated to us by Bence Horv\'ath and Niels Laustsen, who have kindly allowed us to include the details here. Recall that if $\Omega_1, \Omega_2$ are compact Hausdorff spaces, then any continuous map $\theta: \Omega_1 \to \Omega_2$ induces a norm one unital algebra homomorphism $\theta^\circ: C(\Omega_2) \to C(\Omega_1)$ by $\theta^\circ (f) = f \circ \theta$. 

For a locally compact Hausdorff space $K$, we define $\alpha K = K \cup \{ \infty_K \}$ to be the one-point compactification of $K$. It is convenient to use the following alternative description of the Banach algebra of continuous functions on $K$ which vanish at infinity, 
\begin{align*}
    C_0(K) &= \{ f|_K : f \in C(\alpha K),\ f(\infty_K) = 0 \}.
\end{align*}

In the case where $\Omega_1 = \alpha K$ and $\Omega_2 = \alpha L$ are the one-point compactifications of some locally compact Hausdorff spaces $K$ and $L$ and $\theta(\infty_K) = \infty_L$, $\theta^\circ$ maps $C_0(L)$ into $C_0(K)$, so $\theta^\circ$ restricts to a map between these algebras, denoted $\theta^\circ_0$.

\begin{lemma}\label{lmm: niels-and-bence}
    Let $L$ be an open subset of a locally compact Hausdorff space $K$. Then the map $\tau \colon \alpha K \to \alpha L$ given by
    \begin{equation*}
        \tau(k) = \begin{cases}
                k & \text{if } k \in L \\
                \infty_L & \text{otherwise}
            \end{cases}
    \end{equation*}
    is a continuous surjection, and hence the induced map $\tau^\circ_0 \colon C_0(L) \to C_0(K)$ is an isometric algebra homomorphism.
\end{lemma}
\begin{proof}
    It is clear that $\tau$ is surjective. We check the continuity of $\tau$ directly from the definition 
    of the topologies of $\alpha L$ and $\alpha K$. First, suppose that $G \subseteq L$ is open. Then we have 
    $\tau^{-1}(G) = G$, which is an open subset of $K$ because $L \subseteq K$ is open, and therefore 
    $\tau^{-1}(G)$ is open in $\alpha K$. Second, suppose that $G = \alpha L \setminus C$, where $C \subseteq L$ is compact. 
    Then we have $\tau^{-1}(G) = \alpha K \setminus C$, which is open in $\alpha K$ because compactness is independent 
    of the superspace, so $C$ is compact as a subspace of $K$.
\end{proof}

\begin{rem}
    Lemma \ref{lmm: niels-and-bence} does not imply that $\alpha L$ is a retract of $\alpha K$. The reason is that $\alpha L$ may not be a topological subspace of $\alpha K$ because the map $\alpha L \to \alpha K$ which is the identity on $L$ and maps $\infty_L$ to $\infty_K$ may be discontinuous at $\infty_L$. This explains why we distinguish the points at infinity by adding the subscripts $K$ and $L$ to them.
\end{rem}

Using the previous result, we can give the promised canonical inclusion of the spaces $C_0(K_\mathcal{A}) \xhookrightarrow{} C_0(K_\mathcal{B})$ for almost disjoint families $\mathcal{A} \subseteq \mathcal{B} \subseteq \inset{\mathbb{N}}$. 

\begin{lemma}\label{lmm: coupling-almost-disjoint-families1}
    Let $\mathcal{A} \subseteq \mathcal{B} \subseteq \inset{\mathbb{N}}$ be almost disjoint families. Then $K_\mathcal{A}$ is an open subset of $K_\mathcal{B}$.
\end{lemma}
\begin{proof}
    It is clear that if we identify $K_\mathcal{A}$ as a topological subspace of $K_\mathcal{B}$ via the natural set inclusion, then the subspace topology that $K_\mathcal{A}$ inherits from $K_\mathcal{B}$ coincides with the topology given by the almost disjoint family $\mathcal{A}$. Moreover, since $K_\mathcal{B} \setminus K_\mathcal{A} = \{y_B: B \in \mathcal{B} \setminus \mathcal{A}\}$ is closed, $K_\mathcal{A}$ is an open subset of $K_\mathcal{B}$. 
\end{proof}

\begin{rem}\label{rem: cannonical-embedding}
    By Lemmas \ref{lmm: niels-and-bence} and \ref{lmm: coupling-almost-disjoint-families1} the induced map $\tau_0^\circ: C_0(K_\mathcal{A}) \to C_0(K_\mathcal{B})$ is an isometric algebra homomorphism. We refer to it as the canonical embedding and denote it with the symbol $\iota$. Suppose that $\mathcal{A} \subseteq \mathcal{B} \subseteq \inset{\mathbb{N}}$ are almost disjoint families. To avoid confusion, we will denote the isolated points of $K_\mathcal{B}$ by $\{z_n : n \in \mathbb{N}\}$; that is $K_\mathcal{A} = \{x_n: n \in \mathbb{N}\} \cup \{y_A: A \in \mathcal{A}  \}$ and $K_\mathcal{B} = \{z_n: n \in \mathbb{N}\} \cup \{y_B: B \in \mathcal{B}  \}$. Then $\iota: C_0(K_\mathcal{A}) \to C_0(K_\mathcal{B})$ can be explicitly described, for any $f \in C_0(K_{\mathcal{A}})$, by
\begin{equation*}
    (\iota f)(z_n) = f(x_n), \qquad
    (\iota f)(y_A) =
    \begin{cases}
    f(y_A), & \text{if } A \in \mathcal{A}, \\
    0, & \text{if } A \in \mathcal{B} \setminus \mathcal{A}.
    \end{cases}
\end{equation*}

Further, observe that if we have nested almost disjoint families $\mathcal{A}_r \subseteq \mathcal{A}_s \subseteq \mathcal{A}_t \subseteq \inset{\mathbb{N}}$ then the corresponding canonical embeddings satisfy $\iota_{s,t} \iota_{r,s} = \iota_{r,t}$.
\end{rem}

\begin{corollary}\label{cor: coupling-almost-disjoint-families2}
    Let $\mathcal{A} \subseteq \inset{\mathbb{N}}$ be an almost disjoint family and $B \in \inset{\mathbb{N}}$ such that $B \not \in \mathcal{A}$ and $\{B\} \cup \mathcal{A}$ is an almost disjoint family. Then for any $f \in C_0(K_\mathcal{A})$ we have $\lim_{k \in B} f(x_k) = 0$.
\end{corollary}
\begin{proof}
    Let $\iota: C_0(K_\mathcal{A}) \to C_0(K_{\mathcal{A} \cup \{B\}})$ be the canonical embedding and let $f \in C_0(K_\mathcal{A})$. It follows that
    \begin{equation*}
        0 = (\iota f)(y_B) = \lim_{k \in B} (\iota f)(z_k) =  \lim_{k \in B} f(x_k). \qedhere
    \end{equation*} 
\end{proof}

Before we move on from topological considerations, we would like to highlight an important point. While continuous functions $\theta: \Omega_1 \to \Omega_2$ between compact Hausdorff spaces induce norm one operators $\theta^\circ: C(\Omega_2) \to C(\Omega_1)$, the same is not true for continuous functions $\theta: K \to L$ between locally compact Hausdorff spaces. That is, if one tries to define an induced operator $\theta^\circ: C_0(L) \to C_0(K)$ by $\theta^\circ(f) = f \circ \theta$, then one runs into the issue that while $f \circ \theta: K \to \mathbb{K}$ is a continuous function, it need not vanish at infinity.

In fact, Theorem \ref{th: almost-disjoint-families-few-operators} gives an example of locally compact spaces $K_\mathcal{A} \subsetneq K_\mathcal{B}$ so that whenever $\theta: K_\mathcal{A} \to K_\mathcal{B}$ is continuous and $\theta^\circ$ maps functions that vanish at infinity to functions that vanish at infinity, then $\theta$ has countable range. Details of how this works are similar to those in \cite[Proposition 36]{koszmider2021banach}.

\subsection{An overview of matrices acting on almost disjoint families.}\label{subsec: matrices} 

We give a concise account of how matrices act on $C_0(K_{\mathcal{A}})$, in a sense that will be made precise below. We limit ourselves to the essential concepts needed to make this paper self-contained, slightly modifying the definitions and results from Koszmider and Laustsen \cite{koszmider2021banach} to handle operators between different $C_0(K_{\mathcal{A}})$-spaces. For a more detailed discussion, we refer the reader to their work.

\begin{definition}
    Let $M = (m_{k,n})_{k,n \in \mathbb{N}}$ be a matrix with entries in $\mathbb{K}$. We define
    \begin{itemize}
        \item $\mathbb{M} = \{ M : \norm{M} < \infty \}$ where $\norm{M} = \sup \{ \sum_{n \in \mathbb{N}} |m_{k,n}|: k \in \mathbb{N} \}$.
        \item For $f \in \ell_\infty$ and $M \in \mathbb{M}$ we define $Mf \in \ell_\infty$ by $Mf(k) = \sum_{n \in \mathbb{N}} m_{k,n}f(n)$ for $k \in \mathbb{N}$.
        \item We let $I$ denote the $\mathbb{N} \times \mathbb{N}$ identity matrix and for $A \subseteq \mathbb{N}$ we denote by $1_A \in \ell_\infty$ the indicator function of the set $A$.
    \end{itemize}
\end{definition}

We now define the natural analogues of admission, rejection, and undermining in the context of different domain and codomain, see \cite[Definition 14]{koszmider2021banach}.

\begin{definition}
Suppose that $\mathcal{A}, \mathcal{B} \subseteq \inset{\mathbb{N}}$ are almost disjoint families, $C,D \in \inset{\mathbb{N}}$, and $M \in \mathbb{M}$. We say that:
\begin{itemize}
    \item $\left(\mathcal{A}, \mathcal{B} \right)$ \emph{admits} $M$ if $\lim_{k \in B} (M1_A)(k) = 0$ for every $A \in \mathcal{A}$, $B \in \mathcal{B}$. If $\mathcal{A} = \mathcal{B}$, we will say that $\mathcal{A}$ admits $M$.
    \item $(C,D)$ \emph{rejects} $M$ if $\lim_{k \in D} (M1_C)(k)$ does not exist. If $C = D$, we will say that $C$ rejects $M$.
    \item $C$ \emph{undermines} $M$ if there is $n \in \mathbb{N}$ such that $((M 1_{\{n\}})(k))_{k \in C}$ does not converge to $0$.
\end{itemize}
\end{definition}

Similarly, we recall the concept of compact matrices and introduce the corresponding definition for acceptance, see \cite[Definition 15]{koszmider2021banach}.

\begin{definition}
    Suppose that $M \in \mathbb{M}$ and that $A, B \in \inset{\mathbb{N}}.$
    \begin{itemize}
        \item For $j \in \mathbb{N}$, we define $M_j = (m'_{k,n})_{k,n \in \mathbb{N}}$ where $m'_{k,n} = 0$ if $n \leq j$ and $m'_{k,n} = m_{k,n}$ otherwise.
        \item We say that $M$ is a \emph{compact matrix} if $\lim_{j \in \mathbb{N}} \norm{M_j} = 0$.
        \item We define $M^A_B = (m'_{k,n})_{k,n \in \mathbb{N}}$, where $m'_{k,n} = m_{k,n}$ if $(k,n) \in B \times A$, and $m'_{k,n} = 0$ otherwise.
        \item We say that $(A,B)$ \emph{accepts} $M$ if $M^A_B$ is a compact matrix.
        \item Finally, for families $\mathcal{A}, \mathcal{B} \subseteq \inset{\mathbb{N}}$ we say that $(\mathcal{A}, \mathcal{B})$ \emph{accepts} $M$ if $(A \cup A', B\cup B')$ accepts $M$ for every $A, A' \in \mathcal{A}$, $B, B' \in \mathcal{B}$. If $\mathcal{A} = \mathcal{B}$, we say that $\mathcal{A}$ accepts $M$.
    \end{itemize}
\end{definition}

\begin{rem}
   We note that our formulation of $\mathcal{A}$ accepting $M$ is formally stronger than that of Koszmider and Laustsen. However, using Monotonicity \cite[Lemma 28]{koszmider2021banach}, the reader can readily verify that the two formulations are in fact equivalent.
\end{rem}

Lastly, we discuss the concepts of Admission and Monotonicity, see \cite[Lemma 28]{koszmider2021banach}. The proof is essentially the same as in the case of operators whose domain and codomain coincide. We include it for the reader's convenience.

\begin{lemma}\label{lmm: admission-monotonicity}
    Let $M \in \mathbb{M}$ and $\mathcal{A}, \mathcal{B} \subseteq \inset{\mathbb{N}}$. Suppose that $(\mathcal{A}, \mathcal{B})$ accepts $M$.
    \begin{enumerate}[label=(\arabic*)]
        \item (Admission) If no $B \in \mathcal{B}$ undermines $M$, then $(\mathcal{A}, \mathcal{B})$ admits $M$.
        \item (Monotonicity) Suppose that $C \in \inset{\mathbb{N}}$ is contained in a finite union of elements of $\mathcal{A}$ and $D \in \inset{\mathbb{N}} $ in a finite union of elements of $\mathcal{B}$. Then $(C,D)$ accepts $M$.
    \end{enumerate}
\end{lemma}
\begin{proof}
    (Admission) Fix $A \in \mathcal{A}$, $B \in \mathcal{B}$. For every $f \in \ell_\infty$ and $j \in \mathbb{N}$, we have
    \begin{equation*}
        (M^{A}_{B} - (M^{A}_{B})_j)f = \sum_{n \leq j} f(n) M^{A}_{B} 1_{\{n\}},
    \end{equation*}
    which belongs to $c_0$ since the hypothesis that $B$ does not undermine $M$ implies that $\lim_{k \in B} (M1_{\{n\}})(k) = 0$ for every $n \in \mathbb{N}$. Also
    \begin{equation*}
        M^A_B f = (M^A_B - (M^A_B)_j)f + (M^A_B)_j f,
    \end{equation*}
    where $\norm{(M^A_B)_j f}_\infty \leq \norm{(M^A_B)_j} \norm{f}_\infty \to 0$ as $j \to \infty$ by the compactness of $M^A_B$. It follows that $M^A_B f$ can be approximated by elements of $c_0$. Therefore $M^A_B f \in c_0$,  since $c_0$ is closed in $\ell_\infty$. Taking $f = 1_A$ we have
    \begin{equation*}
        0 = \lim_{k \in \mathbb{N}}(M^A_B 1_A)(k) = \lim_{k \in B} (M1_A)(k).
    \end{equation*}
    Since $A$ and $B$ were arbitrary, the result follows.
    
    (Monotonicity) Let $A_1, \dots, A_p \in \mathcal{A}$ and $B_1, \dots, B_q \in \mathcal{B}$ such that $C \subseteq A_1 \cup \dots \cup A_p$ and $D \subseteq B_1 \cup \dots \cup B_q$.

    Since $(\mathcal{A}, \mathcal{B})$ accepts $M$, the matrices $M^{A_i}_{B_l}$ are compact for every $1 \leq i \leq p, 1 \leq l \leq q$ and thus we have
    \begin{equation*}
        \lim_{j \in \mathbb{N}} \left(\sup_{k \in B_l} \sum_{n \in A_i, n > j} |m_{k,n}| \right) = 0.
    \end{equation*}
    It follows that
    \begin{equation*}
         \lim_{j \in \mathbb{N}} \left(\sup_{k \in D} \sum_{n \in C, n > j} |m_{k,n}| \right) = 0,
    \end{equation*}
    so that $M^C_D$ is compact.
\end{proof}

\subsection{Operators and their matrix representation.}\label{sub-sec: operators-and-matrices}

We explore how operators $T: C_0(K_\mathcal{A}) \to C_0(K_\mathcal{B})$ can be related to the matrix theory described in Subsection \ref{subsec: matrices}. We introduce the following analogue of \cite[Definition 18]{koszmider2021banach}.

\begin{definition}
    Let $K, L$ be scattered, locally compact Hausdorff spaces and $T: C_0(K) \to C_0(L)$ be an operator. Then the \emph{reduced matrix of $T$} is the matrix $M^r_T = (m_{z,x})_{z \in L \setminus L', x \in K \setminus K'}$ given by $m_{z,x} = T^*\delta_{z}(\{x\})$ for $x \in K \setminus K'$, $z \in L \setminus L'$, where recall that $K \setminus K'$ and $L \setminus L'$ denote the set of isolated points of $K$ and $L$ respectively.
\end{definition}

Just as a continuous function is determined by its values on a dense set, an operator is characterised by its reduced matrix, in a sense that will be made precise in the following lemma, which is the counterpart of \cite[Lemma 19]{koszmider2021banach} allowing for different domains and codomains.

\begin{lemma}\label{lmm: dense-countable}
    Let $K, L$ be infinite, scattered, locally compact Hausdorff spaces, and let $F = K \setminus K'$ and $G = L \setminus L'$ be the sets of isolated points. Suppose $T: C_0(K) \to C_0(L)$ is an operator. Then $K'$ contains a set $E$ of cardinality not bigger than the cardinality of $G$ such that
    \begin{equation*}
        T(f)|G = M^r_T(f|F)
    \end{equation*}
    whenever $f \in C_0(K)$ and $f(y) = 0$ for every $y \in E$.
\end{lemma}

The previous lemma allows us to analyse $T$ by focusing our attention on its reduced matrix $M^r_T$.

\begin{lemma}\label{lmm: all-combined}
    Suppose that $\mathcal{A}, \mathcal{B} \subseteq \inset{\mathbb{N}}$ are almost disjoint families and that $T: C_0(K_\mathcal{A}) \to C_0(K_\mathcal{B})$ is an operator.
    \begin{enumerate}[label = (\alph*), ref = (\alph*)]
        \item \label{it: operator-to-matrix} There exists $\mathcal{A}' \subseteq \mathcal{A}$ countable such that for all $A \in \mathcal{A} \setminus \mathcal{A}'$ and every $k \in \mathbb{N}$ we have
     \begin{equation*}
         T1_{U(A)}(z_k) = (M^r_T 1_A)(k).
     \end{equation*}
     Consequently, no $(A, B) \in (\mathcal{A \setminus A'}) \times \mathcal{B}$ rejects $M^r_T$.
     \item \label{it: admission-countable-range} If there exist countable subsets $\mathcal{A}'' \subseteq \mathcal{A}$, $\mathcal{B}'' \subseteq \mathcal{B}$ such that $(\mathcal{A} \setminus \mathcal{A}'', \mathcal{B} \setminus \mathcal{B}'')$ admits $M^r_T$, then $T$ factors through $c_0$. 
    \end{enumerate}
\end{lemma}
\begin{proof}
    \ref{it: operator-to-matrix}  Since $G = K_\mathcal{B} \setminus K'_\mathcal{B} = \{z_n: n \in \mathbb{N}\}$ is countable, by Lemma \ref{lmm: dense-countable} applied to $K = K_\mathcal{A}$, $L = K_\mathcal{B}$ there is a countable set $E \subseteq K' = \{y_A: A \in \mathcal{A}\}$ such that $T(f)|G = M^r_T (f|F)$ whenever $f \in C_0(K_\mathcal{A})$ satisfies $f(y) = 0$ for every $y \in E$. Letting $\mathcal{A}' = \{A \in \mathcal{A}: y_A \in E\}$ the result follows. 

    It  follows that for any $A \in \mathcal{A} \setminus \mathcal{A}'$ and any $B \in \mathcal{B}$ we have
    \begin{equation*}
        T(1_{U(A)})(y_B) = \lim_{k \in B} T(1_{U(A)})(z_k) = \lim_{k \in B} (M^r_T 1_A)(k), 
    \end{equation*}
    and thus $(A,B)$ does not reject $M^r_T$.

    \ref{it: admission-countable-range} The hypothesis that $(\mathcal{A} \setminus \mathcal{A}'', \mathcal{B} \setminus \mathcal{B}'')$ admits $M^r_T$ implies that
    \begin{equation*}
        0 = \lim_{k \in B} (M^r_T 1_A)(k) = \lim_{k \in B} (T 1_{U(A)})(z_k) = T1_{U(A)}(y_B)
    \end{equation*}
    for every $A \in \mathcal{A} \setminus (\mathcal{A'} \cup \mathcal{A}'')$ and $B \in \mathcal{B} \setminus \mathcal{B}''$. Since the characteristic functions $1_{\{x_n\}}$ and $1_{U(A)}$ span a dense set of $C_0(K_\mathcal{A})$ it follows that
    \begin{equation*}
        s(T[C_0(K_\mathcal{A})]) \subseteq \mathcal{B}'' \cup \left(\bigcup_{n \in \mathbb{N}}s(T1_{\{x_n\}}) \right) \cup \left(\bigcup_{A \in \mathcal{A}' \cup \mathcal{A}''} s(T1_{U(A)}) \right).
    \end{equation*}
    For each $f \in C_0(K_\mathcal{A})$ the set $s(Tf)$ is countable, thus it follows that $s(T[C_0(K_\mathcal{A})])$ is countable. By Lemma \ref{lmm: s-countable} $T[C_0(K_\mathcal{A})]$ is separable and therefore Lemma \ref{lmm: countable-to-fact-c0} gives that $T$ factors through $c_0$.
\end{proof}

\bigskip
\section{Proof of Theorem \ref{th: almost-disjoint-families-few-operators}} \label{sec: proof-main-theorem}

For the proof, we will need the following two key results from the construction of Koszmider and Laustsen.

\begin{lemma}[{\cite[Lemma 30]{koszmider2021banach}}]\label{lmm: almost-disjoint-borel-support}
There is an uncountable, almost disjoint family $\mathcal{A} \subseteq \inset{\mathbb{N}}$ such that the set
    \begin{equation*}
        \{1_A : A \in \mathcal{A}\} \subseteq 2^{\mathbb{N}}
    \end{equation*}
is closed with respect to the product topology.
\end{lemma}

\begin{lemma}[A dichotomy for acceptance and rejection {\cite[Lemma 35]{koszmider2021banach}}]\label{lmm: dichotomy-acceptance-and-rejection}
Suppose that $M \in \mathbb{M}$ and $\mathcal{A} \subseteq \inset{\mathbb{N}}$ is an uncountable, almost disjoint family such that $\{1_A : A \in \mathcal{A}\}$ is a Borel subset of $2^{\mathbb{N}}$. Then one of the following holds:
\begin{enumerate}[label = (\alph*), ref = (\alph*)]
    \item either $\mathcal{A} \setminus \mathcal{A}'$ accepts $M - \lambda I$ for some countable $\mathcal{A}' \subseteq \mathcal{A}$ and $\lambda \in \mathbb{K}$,
    \item \label{it: alternative-b} or there are pairwise disjoint sets $\{A_\xi, A'_\xi\} \subseteq \mathcal{A}$ and an infinite subset $B_\xi \subseteq A_\xi \cup A'_\xi$ such that $B_\xi$ rejects $M$ for every $\xi < \mathfrak{c}$.
\end{enumerate}
\end{lemma}

\begin{proof}[Proof of Theorem \ref{th: almost-disjoint-families-few-operators}]
    We follow closely the proof of \cite[Theorem 2]{koszmider2021banach}, doing modifications when necessary.
    \begin{enumerate}[label=(\arabic*), ref = (\arabic*)]
        \item \label{step:1} Using Lemma \ref{lmm: almost-disjoint-borel-support}, we fix an uncountable almost disjoint family $\mathcal{A} \subseteq \inset{\mathbb{N}}$ such that $\{1_A : A \in \mathcal{A}\}$ is a Borel subset of $2^{\mathbb{N}}$.
        \item \label{step:2} Let $\mathbb{M}'$ be the set of all matrices $M \in \mathbb{M}$ such that there is no countable $\mathcal{A}_M \subseteq \mathcal{A}$ and no $\lambda_M \in \mathbb{K}$ such that $\mathcal{A} \setminus \mathcal{A}_M$ accepts $M - \lambda_M I$.
        
        If $\mathbb{M}' = \emptyset$, define $\mathcal{B} = \emptyset$ and $\mathcal{C} = \mathcal{A}$. Move directly to step \ref{it: four}. Otherwise let $\{(M_\xi, \varepsilon_\xi): \xi < \mathfrak{c}\}$ be an enumeration of $\mathbb{M}' \times \{0,1\}$ with each element repeated continuum many times. Observe that for $\mathcal{A}$ and for each $M_\xi$ we are in alternative \ref{it: alternative-b} of the dichotomy for acceptance and rejection, Lemma \ref{lmm: dichotomy-acceptance-and-rejection}.
        \item  By transfinite recursion on $\xi < \mathfrak{c}$, construct $A_\xi, A'_\xi \in \mathcal{A}$, $B_\xi\in \inset{\mathbb{N}}$ such that:
        \begin{enumerate}
            \item $\{A_\eta, A'_\eta \} \cap \{A_\xi, A'_\xi \} = \emptyset$ for all $\eta < \xi$,
            \item $B_\xi \subseteq A_\xi \cup A'_\xi$,
            \item $B_\xi$ rejects $M_\xi$.
        \end{enumerate}
        Let $\mathcal{B} = \{B_\xi: \xi < \mathfrak{c}, \varepsilon_\xi = 0\}$ and $\mathcal{C} = \{B_\xi: \xi < \mathfrak{c}, \varepsilon_\xi = 1\}$. Observe that since $\mathcal{A}$ is an almost disjoint family, so is $\mathcal{B} \cup \mathcal{C}$.

    \item \label{it: four} Consider a partition $\mathcal{C}=\bigsqcup_{\alpha<\mathfrak c}\mathcal{C}_\alpha$ with $|\mathcal{C}_\alpha|=\mathfrak c$ for each $\alpha<\mathfrak c$, and an order embedding $\mathcal P\to\wp(\mathfrak c)$, $t\mapsto\kappa_t$. For each $t\in\mathcal P$, define
    \begin{equation*}
        \mathcal{B}_t=\mathcal{B}\cup\left(\bigcup_{\alpha\in\kappa_t}\mathcal{C}_\alpha\right).
    \end{equation*}
    Observe that, for each $t\in\mathcal P$, $\mathcal{B}_t\subseteq[\mathbb N]^\omega$ is an almost disjoint family because $\mathcal{B}\cup\mathcal{C}$ is. Furthermore, for each $s,t\in\mathcal P$, if $s\leq t$, then $\mathcal{B}_s\subseteq\mathcal{B}_t$, since $\kappa_s\subseteq\kappa_t$. On the other hand, if $s\not\leq t$, then $|\mathcal{B}_s\setminus\mathcal{B}_t|=\mathfrak c$. Indeed, $s\not\leq t$ implies $\kappa_s\not\subseteq\kappa_t$, so we may choose $\alpha\in\kappa_s\setminus\kappa_t$. Then
    \begin{equation*}
        \mathcal{C}_\alpha\subseteq\mathcal{B}_s\setminus\mathcal{B}_t.
    \end{equation*}
    Since $|\mathcal{C}_\alpha|=\mathfrak c$, it follows that $|\mathcal{B}_s\setminus\mathcal{B}_t|=\mathfrak c$.
    
    \item \label{it: five} We check that the collection $(\mathcal{A}_t)_{t\in\mathcal P}=(\mathcal{B}_t)_{t\in\mathcal P}$ satisfies the theorem. Thus, fix $s,t\in\mathcal P$ and an operator $T:C_0(K_{\mathcal{B}_s})\to C_0(K_{\mathcal{B}_t})$.

        We claim that the sets $\mathcal{B}'_t = \{B \in \mathcal{B}_t: B \text{ undermines } M^r_T \}$ and $\mathcal{B}'_s = \{B \in \mathcal{B}_s: B \text{ undermines } M^r_T \}$ are countable. For this, it is enough to show that for each fixed $n \in \mathbb{N}$ the sets
        \begin{equation*}
            V(n) = \{B \in \mathcal{B}_s : ((M^r_T 1_{\{n\}})(k))_{k \in B} \text{ does not converge to } 0 \}
        \end{equation*}
        and
        \begin{equation*}
             W(n) = \{B \in \mathcal{B}_t : ((M^r_T 1_{\{n\}})(k))_{k \in B} \text{ does not converge to } 0 \}
        \end{equation*}
        are countable. Observe that $(M^r_T 1_{\{n\}})(k) = T1_{\{x_n\}}(z_k)$ and thus for $B \in \mathcal{B}_t$ we get
        \begin{equation*}
            \lim_{k \in B} (M^r_T 1_{\{n\}}) (k) = \lim_{k \in B} T1_{\{x_n\}}(z_k) = T1_{\{x_n\}}(y_B).
        \end{equation*}
        Since $T1_{\{x_n\}} \in C_0(K_{\mathcal{B}_t})$,
        it follows that $W(n)$ is countable. 
        
        We claim that $V(n) \subseteq W(n)$, and thus $V(n)$ is also countable. Let $B \in V(n)$, we will argue that $B \in \mathcal{B}_t$ and thus $B \in W(n)$. Indeed, notice that whenever $B \in \mathcal{B}_s$ and $B \not\in \mathcal{B}_t$ then $\{B\} \cup \mathcal{B}_t$ is an almost disjoint family. Since $T1_{\{x_n\}} \in C_0(K_{\mathcal{B}_t})$, Corollary \ref{cor: coupling-almost-disjoint-families2} gives
            \begin{equation*}
                0 = \lim_{k \in B} T1_{\{x_n\}}(z_k) = \lim_{k \in B} (M^r_T 1_{\{n\}})(k),
            \end{equation*}
        so $B \not \in V(n)$.
        
        Observe that for each $M \in \mathbb{M}'$, there exist uncountably many $B \in \mathcal{B} \subseteq \mathcal{B}_s \cap \mathcal{B}_t$ such that $B$ rejects $M$. Thus, by Lemma \ref{lmm: all-combined} \ref{it: operator-to-matrix} we conclude that $M^r_T \not \in \mathbb{M}'$. From the choices made in \ref{step:2}, it follows that we can find $\mathcal{A}_{M^r_T} \subseteq \mathcal{A}$ countable and $\lambda_{M^r_T} \in \mathbb{K}$ such that $\mathcal{A} \setminus \mathcal{A}_{M^r_T}$ accepts $M^r_T - \lambda_{M^r_T} I$. For notational convenience, we will write $\mathcal{A}_T$ for $\mathcal{A}_{M^r_T}$ and $\lambda_T$ for $\lambda_{M^r_T}$.
         
        Since all but countably many elements of $\mathcal{B}_s$ and $\mathcal{B}_t$ are in the union of finitely many elements of $\mathcal{A}\setminus \mathcal{A}_T$, the Monotonicity of Lemma \ref{lmm: admission-monotonicity} implies that we can find countable subsets $\mathcal{B}''_s \subseteq \mathcal{B}_s$, $\mathcal{B}''_t \subseteq \mathcal{B}_t$ such that both the family $\mathcal{B}_s \setminus (\mathcal{B}'_s \cup \mathcal{B}''_s)$ and the pair $(\mathcal{B}_s \setminus (\mathcal{B}'_s \cup \mathcal{B}''_s), \mathcal{B}_t \setminus (\mathcal{B}'_t \cup \mathcal{B}''_t))$ accept $M^r_T - \lambda_T I$. 

        Moreover, since no element of $(\mathcal{B}_s \setminus \mathcal{B}'_s) \cup (\mathcal{B}_t \setminus \mathcal{B}'_t)$ undermines $M^r_T$, they also do not undermine $M^r_T - \lambda_T I$. In particular, by Admission of Lemma \ref{lmm: admission-monotonicity}, it follows that both the family $\mathcal{B}_s \setminus (\mathcal{B}'_s \cup \mathcal{B}''_s)$ and the pair $(\mathcal{B}_s \setminus (\mathcal{B}'_s \cup \mathcal{B}''_s), \mathcal{B}_t \setminus (\mathcal{B}'_t \cup \mathcal{B}''_t))$ admit $M^r_T - \lambda_T I$. 

        To finish the proof, we distinguish two cases: $s \leq t$ and $s \not \leq t$.

        \begin{enumerate}
            \item First, assume that $s \leq t$. By construction $\mathcal{B}_s \subseteq \mathcal{B}_t$ and thus, as in Remark \ref{rem: cannonical-embedding}, we have a canonical embedding $\iota_{s,t}: C_0(K_{\mathcal{B}_s}) \to C_0(K_{\mathcal{B}_t})$, and naturally $M_{\iota_{s,t}}^{r} = I$. It follows that the operator $S = T - \lambda_T \iota_{s,t}$ satisfies the conditions of Lemma \ref{lmm: all-combined} \ref{it: admission-countable-range} and therefore factors through $c_0$.
            \item Assume that $s \not \leq t$. According to Lemma \ref{lmm: all-combined} \ref{it: operator-to-matrix} we can find a countable subset $\mathcal{B}'''_{s} \subseteq \mathcal{B}_{s}$ such that for any $B \in \mathcal{B}_{s} \setminus \mathcal{B}'''_{s}$ we have $T(1_{U(B)})(z_k) = (M^r_T 1_B)(k)$ for all $k \in \mathbb{N}$.

            Observe that $\mathcal{B}'_s$, $\mathcal{B}''_s$ and $\mathcal{B}'''_s$ are countable, while the set $\mathcal{B}_s \setminus \mathcal{B}_t$ is uncountable. Therefore, we can find $B \in \mathcal{B}_s \setminus (\mathcal{B}_t \cup \mathcal{B}'_s \cup \mathcal{B}''_s \cup \mathcal{B}'''_s)$.

            Since $M^r_T - \lambda_T I$ admits $\mathcal{B}_s \setminus (\mathcal{B}'_s \cup \mathcal{B}''_s)$, we have
        \begin{equation*}
            0 = \lim_{k \in B} ((M^r_T - \lambda_T I)1_B) (k),
        \end{equation*}
        or equivalently
        \begin{equation*}
            \lambda_T = \lim_{k \in B} (M^r_T 1_B)(k). 
        \end{equation*}
        Since $B \not \in \mathcal{B}'''_s$, it follows that
        \begin{equation*}
            \lim_{k \in B} (T1_{U(B)})(z_k) = \lim_{k \in B} (M^r_T 1_B)(k )= \lambda_T.
        \end{equation*}
        Since $B \not \in \mathcal{B}_t$, $\{B\} \cup \mathcal{B}_t$ is an almost disjoint family and $T1_{U(B)} \in C_0(K_{\mathcal{B}_t})$ an application of Corollary \ref{cor: coupling-almost-disjoint-families2} gives
        \begin{equation*}
            0 = \lim_{k \in B} (T1_{U(B)})(z_k) = \lambda_T.
        \end{equation*}
        Finally, since $(\mathcal{B}_s \setminus (\mathcal{B}'_s \cup \mathcal{B}''_s), \mathcal{B}_t \setminus (\mathcal{B}'_t \cup \mathcal{B}''_t))$ admits $M^r_T - \lambda_T I = M^r_T$, it follows from Lemma \ref{lmm: all-combined} \ref{it: admission-countable-range} that $T$ factors through $c_0$. This finishes the proof.\qedhere
        \end{enumerate}
    \end{enumerate}
\end{proof}
\bigskip
\section{Quotient algebras and operator ideals}\label{sec: operator-ideals}

We explore some consequences of Theorem \ref{th: almost-disjoint-families-few-operators} for the construction of $C(K)$-spaces whose closed operator ideals can be classified. For convenience, we present our result in slightly more general terms. We will need the following definition.

\begin{definition}\label{def: family-ordered-few}
    Let $\mathcal{P}$ be a partially ordered set. We say that a family $(X_t)_{t \in \mathcal{P}}$ of Banach spaces \emph{has few operators relative to $\mathcal{P}$} if the following are satisfied for $s,t \in \mathcal{P}$:
    \begin{enumerate}[label = (\alph*), ref = (\alph*)]
        \item $X_t$ is not separable.
        \item If $s \leq t$ then there exists an isometric embedding $\iota_{s,t}: X_s \to X_t$ so that every operator $T: X_s \to X_t$ is of the form $T = \lambda \iota_{s,t} + S$, where $\lambda \in \mathbb{K}$ and $S$ has separable range. Furthermore, for every $r \leq s \leq t$ we have $\iota_{s,t} \iota_{r,s} = \iota_{r,t}$.
        \item \label{it: condc} If $s \not \leq t$, then every operator $T: X_s \to X_t$ has separable range.
    \end{enumerate}
\end{definition}

\begin{rem}
    We think of the spaces $(X_t)_{t \in \mathcal{P}}$ as being nested relative to the poset $\mathcal{P}$. Accordingly, whenever we write $X_s \subseteq X_t$ for $s \leq t$, we implicitly identify $X_s$ with its image $\iota_{s,t}(X_s) \subseteq X_t$.
\end{rem}

\begin{rem}
    In case \ref{it: condc}, that is, whenever $s \nleq t$, we adopt the convention of defining $\iota_{s,t} = 0$ to maintain consistency in the notation.
\end{rem}

According to Theorem \ref{th: almost-disjoint-families-few-operators}, for any poset $\mathcal{P}$ of size at most the continuum, we can always find a family of non-separable $C(K)$-spaces which has few operators relative to $\mathcal{P}$. We recall that a poset $\mathcal{P}$ is said to be \emph{locally finite} if for any $s,t \in \mathcal{P}$ with $s \leq t$, the interval $[s,t] = \{r \in \mathcal{P}: s \leq r \leq t \}$ is finite. 

\begin{definition}\label{def: inc-algebra}
    Let $(\mathcal{P}, \leq)$ be a locally finite partially ordered set. The \emph{incidence algebra} of $\mathcal{P}$ is the algebra
    \begin{equation*}
        \incidence(\mathcal{P}) = \{f\colon \mathcal{P} \times \mathcal{P} \to \mathbb{K} : f(s,t) = 0 \text{ for } s,t \in \mathcal{P} \text{ with } s \not\leq t  \},
    \end{equation*}
    with the sum of functions and scalar multiplication defined pointwise. The product $h = f \star g$ is defined via convolution
    \begin{equation*}
        h(s,t) = \sum_{s \leq r \leq t} f(s,r)g(r,t).
    \end{equation*}
\end{definition}

\begin{rem}
    For technical reasons, we will work with the incidence algebra on the poset with the opposite order $\mathcal{P}^{\opposite}$. Observe that this algebra can be explicitly described in terms of the original order as
    \begin{equation*}
        \incidence(\mathcal{P}^{\opposite}) = \{f\colon \mathcal{P} \times \mathcal{P} \to \mathbb{K} : f(s,t) = 0 \text{ for } s,t \in \mathcal{P} \text{ with } t \not\leq s\},
    \end{equation*}
    and the product becomes
    \begin{equation*}
        ( f \star g)(s,t) = \sum_{t \leq r \leq s} f(s,r)g(r,t).
    \end{equation*}
\end{rem}

Incidence algebras were originally introduced by Rota \cite{rota1964foundations} to provide a general framework for inversion-type formulas and have since found multiple applications. Examples of incidence algebras include, among others, the algebra of $n \times n$ upper-triangular matrices $\mathbb{T}_n$, which is induced by the linearly ordered set of $n$ elements.

We would like to point out that incidence algebras naturally carry a topology, given by pointwise convergence of functions. For $\mathcal{P}$ finite, this topology coincides with the one induced by any norm. In this case, we will regard elements in $\incidence(\mathcal{P})$ as operators on $(\mathbb{K}^{\mathcal{P}}, \norm{\cdot}_\infty)$ and equip $\incidence(\mathcal{P})$ with the operator norm. Specifically, $f \in \incidence(\mathcal{P})$ acts on an element $v \in \mathbb{K}^{\mathcal{P}}$ by
\begin{equation*}
    (f \cdot v)(s) = \sum_{t \in \mathcal{P}} f(s,t) v(t).
\end{equation*}
Of special interest in our case is the fact that, with this topology understood, there is a procedure to determine the closed ideals of incidence algebras; see \cite[Section 3]{doubilet1972foundations}.

\subsection{Finite-dimensional incidence algebras}  We shall now show how to create incidence algebras of finite posets as quotients of algebras of operators $\mathscr{B}(X)$ on a Banach space $X$. We have the following elementary observation.

\begin{proposition}\label{prop: create-incidence-algebras}
Suppose that $\mathcal{P}$ is a finite partially ordered set, $(X_t)_{t \in \mathcal{P}}$ is a family that has few operators relative to $\mathcal{P}$ and $Y = \bigoplus_{t \in \mathcal{P}} X_t$. Assume moreover that the spaces $X_t$ contain a common
non-separable subspace $Z$, and that $\iota_{t,s}$ restricts
to the identity on $Z$ whenever $t\leq s$.
Then
\begin{equation*}
    \{0\} \longrightarrow \mathscr{X}(Y) \longrightarrow \mathscr{B}(Y) \xlongrightarrow{\theta} \incidence(\mathcal{P}^{\opposite}) \longrightarrow \{0\}
\end{equation*}
is a strongly split exact sequence, the right inverse $\rho$ of $\theta$ is an isometry and $\norm{\theta} = 1$.
\end{proposition}
\begin{proof}
    For each $t \in \mathcal{P}$, let $\pi_t: Y \to X_t$ be the canonical projection onto $X_t$ and $\jmath_t: X_t \to Y$ be the canonical embedding into $Y$. Any operator $T: Y \to Y$ can be expressed as
    \begin{equation*}
        T = \sum_{s,t \in \mathcal{P}} \jmath_s T_{s,t} \pi_t, \hspace{5pt} \text{ where } \hspace{5pt}T_{s,t} = \pi_s T\jmath_t: X_t \to X_s.
    \end{equation*}
    
    Since the family $(X_t)_{t \in \mathcal{P}}$ has few operators relative to $\mathcal{P}$, we have a unique decomposition $T_{s,t} = \lambda^T_{s,t} \iota_{t,s} + S^T_{s,t}$ where $S^T_{s,t}$ has separable range and $\lambda^T_{s,t} = 0$ whenever $t \not \leq s$. Define $S^T = \sum_{s,t \in \mathcal{P}} \jmath_s S^T_{s,t}\pi_t$, which is an operator with separable range. We get
    \begin{equation*}
        T = \sum_{s,t \in \mathcal{P}, t \leq s} \jmath_s (\lambda^T_{s,t} \iota_{t,s}) \pi_t + \sum_{s,t \in \mathcal{P}} \jmath_s S^T_{s,t}\pi_t = \sum_{s,t \in \mathcal{P}, t \leq s} \jmath_s (\lambda^T_{s,t} \iota_{t,s}) \pi_t + S^T.
    \end{equation*}
    Therefore, we can define the map $\theta: \mathscr{B}(Y) \to \incidence(\mathcal{P}^{\opposite})$ by $\theta(T)(s,t) = \lambda_{s,t}^T$. It is clear that $\theta$ is linear and $\ker \theta = \mathscr{X}(Y)$, we need to check that $\theta$ is multiplicative. Let $T, U \in \mathscr{B}(Y)$, we have
    \begin{equation*}
        TU = \left(\sum_{s,q \in \mathcal{P}, q \leq s} \jmath_s (\lambda^T_{s,q} \iota_{q,s}) \pi_q + S^T \right) \left(\sum_{r,t \in \mathcal{P}, t \leq r} \jmath_r (\lambda^U_{r,t} \iota_{t,r}) \pi_t + S^U \right),
    \end{equation*}
    so that if we define the separable range operator $S^{T,U}$ by
    \begin{equation*}
        S^{T,U} = S^T \left(\sum_{r,t \in \mathcal{P}, t \leq r} \jmath_r (\lambda^U_{r,t} \iota_{t,r}) \pi_t \right) + \left(\sum_{s,q \in \mathcal{P}, q \leq s} \jmath_s (\lambda^T_{s,q} \iota_{q,s}) \pi_q \right) S^U + S^T S^U,
    \end{equation*}
    then we have
    \begin{align*}
        TU &= \left(\sum_{s,q \in \mathcal{P}, q \leq s} \jmath_s (\lambda^T_{s,q} \iota_{q,s}) \pi_q \right) \left(\sum_{r,t \in \mathcal{P}, t \leq r} \jmath_r (\lambda^U_{r,t} \iota_{t,r}) \pi_t\right) + S^{T,U} \\
        &= \left(\sum_{s,t \in \mathcal{P}} \jmath_s \left(\sum_{t \leq r \leq s} \lambda^T_{s,r}  \lambda^U_{r,t} \right) \iota_{t,s} \pi_t\right) + S^{T,U} 
        \\&= \left(\sum_{s,t \in \mathcal{P}} \jmath_s (\theta(T) \star \theta(U))(s,t) \iota_{t,s} \pi_t\right) + S^{T,U},
    \end{align*}
    and thus $\theta(TU) = \theta(T) \star \theta(U)$.
    
    It is clear that the map
    \begin{equation*}
    \rho: \incidence(\mathcal{P}^{\opposite}) \longrightarrow \mathscr{B}(Y), \qquad
    \rho(f) = \sum_{s,t \in \mathcal{P}} \jmath_s f(s,t) \iota_{t,s} \pi_t,
    \end{equation*}
    defines an algebra homomorphism which is a right inverse to $\theta$. 
    Since the finite direct sum $Y = \bigoplus_{t \in \mathcal{P}} X_t$ is taken with respect to the $\|\cdot\|_\infty$-norm, elementary arguments using the common non-separable subspace $Z$ and the finiteness of $\mathcal{P}$ show that $\rho$ is an isometry and that $\norm{\theta} = 1$.
\end{proof}

The choice of the ideal of separable-range operators is not essential in the finite-dimensional case. The same argument applies, for instance, to compact operators, provided the corresponding scalar-plus-compact decomposition is available. This was implicitly done, for the two-point linearly ordered set, by Kania and Laustsen \cite{kania2017ideal}, who obtained the upper-triangular $2 \times 2$ matrices as a Calkin algebra. This suggests the following related question.

\begin{question}
Given a finite partially ordered set $\mathcal{P}$, is the incidence algebra $\incidence(\mathcal{P})$ a Calkin algebra?
\end{question}

A positive answer would follow by constructing, for each finite $\mathcal{P}$, Argyros-Haydon-like spaces with few operators relative to $\mathcal{P}$, with compact operators replacing separable-range operators.

\subsection{$c_0$-Incidence algebras.}

We aim to expand the previous result to the case of infinite incidence algebras. For countable partially ordered sets, the definition of incidence algebra requires local finiteness of the poset $\mathcal{P}$ so that the convolution product is well-defined. Unfortunately, even under the assumption of local finiteness of $\mathcal{P}$, the quotient $\mathscr{B}(X)/\mathscr{X}(X)$ will never be one of these algebras. We have the following stronger statement.

\begin{proposition}
    Let $(\mathcal{P}, \leq)$ be a locally finite partially ordered set. Then the incidence algebra $\incidence(\mathcal{P})$ is normable if and only if $\mathcal{P}$ is finite.
\end{proposition}
\begin{proof}
    If $\mathcal{P}$ is finite, then $\incidence(\mathcal{P})$ is a subalgebra of a  finite-dimensional matrix algebra, so any norm on the latter gives a norm on $\incidence(\mathcal{P})$. 
    
    Suppose now that $\mathcal{P}$ is infinite and let $\phi: \mathbb{N} \to \mathcal{P}$ be an injection. Define the functions $f: \mathcal{P} \times \mathcal{P} \to \mathbb{K}$ and $g_n : \mathcal{P} \times \mathcal{P} \to \mathbb{K}$ for each $n \in \mathbb{N}$ by
    \begin{equation*}
    f(s,t) = 
    \begin{cases}
    n & \text{if } s = t = \phi(n) \text{ for some } n, \\
    0 & \text{otherwise},
    \end{cases}
    \qquad
    g_n(s,t) =
    \begin{cases}
    1 & \text{if } s = t = \phi(n), \\
    0 & \text{otherwise}.
    \end{cases}
    \end{equation*} 
    Clearly $f, g_n \in \incidence(\mathcal{P})$, while $f \star g_n = n g_n$. Now the conclusion follows from \cite[Proposition 2.1.3]{dales2001banach}. Indeed, suppose that $\norm{\cdot}$ is a submultiplicative norm on $\incidence(\mathcal{P})$, then
    \begin{equation*}
        n \norm{g_n} = \norm{ng_n} = \norm{f \star g_n} \leq \norm{f} \norm{g_n},
    \end{equation*}
    that is, $n \leq \norm{f}$ for every $n \in \mathbb{N}$, which gives the desired contradiction.
    \end{proof}

This section aims to introduce an alternative to incidence algebras that is better suited to our context and can be realised as quotients of algebras of operators.

In the remainder of this section, we fix a countably infinite partially ordered set $(\mathcal{P}, \leq )$. We will denote by $(e_t)_{t\in \mathcal{P}}$ the usual basis of $c_0(\mathcal{P})$ and $(e^*_t)_{t\in \mathcal{P}}$ the biorthogonal functionals. 

\begin{definition}\label{def: E-inc-algebra}
    The \emph{$c_0$-incidence algebra} $\incidence(\mathcal{P}, c_0)$ is the subspace of $\mathscr{B}(c_0(\mathcal{P}))$ given by operators $T$ satisfying $e^*_s T e_t = 0$ whenever $s \not\leq t$ for $s,t \in \mathcal{P}$.
\end{definition}
\begin{rem}
    We observe that, in principle, one can define $E$-incidence algebras similarly for any Banach space $E$ with an unconditional basis. However, since such algebras are not well-suited to our purposes, we will not pursue this direction.
\end{rem}

We begin by verifying that the above definition yields a Banach algebra.

\begin{proposition}
     $\incidence(\mathcal{P}, c_0)$ is a closed subalgebra of $\mathscr{B}(c_0(\mathcal{P}))$.
\end{proposition}
\begin{proof}
    From the definition, it is immediate that $\incidence(\mathcal{P}, c_0)$ is a closed subspace. Thus, we only need to verify that it is closed under multiplication. For each $r \in \mathcal{P}$, let $e_r\otimes e_r^*$ denote the rank-one projection given by $(e_r\otimes e_r^*)(x)=e_r^*(x)e_r$ for $x \in c_0(\mathcal{P})$. Observe that $I=\sum_{r\in\mathcal P}e_r\otimes e_r^*$, with convergence in the strong operator topology. Let $S, T \in \incidence(\mathcal{P}, c_0)$, and choose $s, t$ such that $s \not\leq t$. Then
    \begin{equation*}
        e_s^* S T e_t = e_s^* S \left(\sum_{r \in \mathcal{P}} e_r \otimes e_r^* \right) T e_t = \sum_{r \in \mathcal{P}} (e_s^* S e_r)(e_r^*Te_t) = 0,
    \end{equation*}
    since whenever $s \not \leq r$ we have $e_s^* S e_r = 0$, while when $s \leq r$ then $r \not \leq t$ since $s \not\leq t$, so that $e_r^*Te_t = 0$.
\end{proof}

To each operator $T \in \mathscr{B}(c_0(\mathcal{P}))$ we can associate the infinite matrix $(\lambda_{s,t})_{s,t \in \mathcal{P}}$ by $\lambda_{s,t} = e^*_s T e_t$. We call this matrix the \emph{scalar matrix of $T$}. Under this identification,  $\incidence(\mathcal{P},c_0)$ is the subalgebra of operators on $c_0(\mathcal{P})$ whose matrix representation has zeros according to the order of $\mathcal{P}$.

\begin{example}
    Let $\mathcal{P} = \mathbb{N}$ with the usual order $\leq$. Then $\incidence(\mathcal{P},c_0)$ is the subalgebra of operators on $c_0$ that have upper-triangular matrix representation. Similarly, if we take $\mathcal{P} = \mathbb{N}$ with the opposite order relation $\geq$, we get the algebra of operators that have lower-triangular matrix representation.
\end{example}

Our goal is to give an analogue of Proposition \ref{prop: create-incidence-algebras} for $c_0$-incidence algebras. For technical reasons, it is useful to assume the existence of a superspace $X_\infty$ that contains all members of the family, as well as a minimal space $X_0$ that is contained in all of them, so that $X_0 \subseteq X_t \subseteq X_\infty$ for all $t \in \mathcal{P}$. These considerations motivate the following definition.

\begin{definition}\label{def: bounded-extension}
    We define the \emph{bounded extension} of $\mathcal{P}$ to be the poset $(\mathcal{P}^*, \leq)$ where we adjoin two formal extra elements $\mathcal{P}^* = \mathcal{P} \cup \{0, \infty\}$ and declare $s \leq t$ whenever $t = \infty$, or $s = 0$ or $s, t \in \mathcal{P}$ and $s \leq t$.
\end{definition}

\begin{definition}
    Let $\Lambda = (\lambda_{s,t})_{s,t \in \mathcal{P}}$ be an infinite scalar matrix. We say that $\Lambda$ has $\ell_1$-\emph{bounded rows} if
        \begin{equation*}
            \sup_{s \in \mathcal{P}} \sum_{t \in \mathcal{P}} |\lambda_{s,t}| < \infty.
        \end{equation*}
    We say that $\Lambda$ has \emph{$c_0$-columns} if for every $t \in \mathcal{P}$ we have $(\lambda_{s,t})_{s \in \mathcal{P}} \in c_0(\mathcal{P})$.
\end{definition}

We recall the following well-known characterization of when an infinite scalar matrix $(\lambda_{s,t})_{s,t \in \mathcal{P}}$ defines an operator in $\mathscr{B}(c_0(\mathcal{P}))$, see for example \cite[Section IV.6, Theorem 6.3]{taylor1980introduction}.

\begin{proposition}\label{prop: when-matrix-is-operator}
    \begin{enumerate}[label = (\alph*), ref = (\alph*)]
        \item \label{it: when-matrix-is-operator-a} Let $T \in \mathscr{B}(c_0(\mathcal{P}))$. Then $(\lambda_{s,t})_{s,t \in \mathcal{P}}$, the scalar matrix of $T$, has $\ell_1$-bounded rows and $c_0$-columns. Moreover, 
        \begin{equation*}
            \norm{T} = \sup_{s \in \mathcal{P}} \sum_{t \in \mathcal{P}} |\lambda_{s,t}|.
        \end{equation*}
        \item \label{it: when-matrix-is-operator-b} Let $\Lambda = (\lambda_{s,t})_{s,t \in \mathcal{P}}$ be an infinite scalar matrix having $\ell_1$-bounded rows and $c_0$-columns. Then there exists an operator $T \in \mathscr{B}(c_0(\mathcal{P}))$ with $\Lambda$ as its scalar matrix.
    \end{enumerate}
\end{proposition}

If $(X_t)_{t \in \mathcal{P}}$ is a family of Banach spaces, we define
\begin{equation*}
    Y = \left( \bigoplus_{t \in \mathcal{P}} X_t\right)_{c_0} = \{(x_t)_{t \in \mathcal{P}}: x_t \in X_t \text{ for } t \in \mathcal{P} \text{ and } (\norm{x_t})_{t \in \mathcal{P}} \in c_0(\mathcal{P})\}.
\end{equation*}
Observe that any operator $T \in \mathscr{B}(Y)$ can be identified with an infinite matrix of operators $(T_{s,t})_{s,t \in \mathcal{P}}$ given by $T_{s,t} = \pi_s T \jmath_t \in \mathscr{B}(X_t, X_s)$, where $\pi_s$ is the canonical projection of $Y$ onto $X_s$ and $\jmath_t$ the canonical inclusion of $X_t$ into $Y$.

\begin{rem}
    As in the case of incidence algebras, for technical reasons, it will be convenient to work with the $c_0$-incidence algebra in the opposite order. In terms of the original order, this means that the subalgebra of operators is defined by the condition $e_s^* T e_t = 0$ whenever $t \nleq s$ for $s, t \in \mathcal{P}$. Note that, since the underlying set is unchanged, we have $c_0(\mathcal{P}) = c_0(\mathcal{P}^{\opposite})$.
\end{rem}

\begin{lemma}\label{lmm: from-scalar-going-up}
   Let $(X_t)_{t \in \mathcal{P}^*}$ be a family of Banach spaces with few operators relative to $\mathcal{P}^*$, and set $Y = \left(\bigoplus_{t \in \mathcal{P}} X_{t} \right)_{c_0}$. Then the map $\rho: \incidence(\mathcal{P}^{\opposite}, c_0) \to  \mathscr{B}(Y)$, sending an operator $T$ with scalar matrix $\Lambda = (\lambda_{s,t})_{s,t \in \mathcal{P}}$ to the operator with matrix $(\lambda_{s,t}\iota_{t,s})_{s,t \in \mathcal{P}}$ defines an isometric homomorphism.
\end{lemma}
\begin{proof}
    Since $T$ is in $\incidence(\mathcal{P}^{\opposite}, c_0)$,  $\lambda_{s,t} \not = 0$ implies that $t \leq s$. Now, whenever $t \leq s$ we have isometries $\iota_{t,s}: X_t \to X_s$ since $(X_t)_{t \in \mathcal{P}}$ has few operators relative to $\mathcal{P}$. This implies that the matrix $(\lambda_{s,t}\iota_{t,s})_{s,t \in \mathcal{P}}$ is well-defined. A straightforward computation, using that $\Lambda$ has $\ell_1$-bounded rows and $c_0$-columns, shows that $(\lambda_{s,t}\iota_{t,s})_{s,t \in \mathcal{P}}$ represents an operator $\rho(T) \in \mathscr{B}(Y)$. That $\rho$ is an isometric homomorphism is elementary, and thus we omit the details.
\end{proof}

We will need the following easy observations.

\begin{lemma}\label{lmm: separable-to-separable}
    Let $X, Y$ be Banach spaces and $S_n: X \to Y$ be an operator with separable range for each $n \in \mathbb{N}$. Suppose that $S_n \to S$ in the strong operator topology. Then $S$ has separable range.
\end{lemma}
\begin{proof}
    Let $Z = \overline{\spn} \{S_n[X]: n \in \mathbb{N} \}$ which is separable. We claim that $S[X] \subseteq Z$, which will finish the proof. Let $y \in S[X]$ and choose $x \in X$ such that $Sx = y$. Since $S_n \to S$ in the strong operator topology, we have $\lim_{n \to \infty} S_n x = y$, and since $S_n x \in Z$ for each $n \in \mathbb{N}$, the result follows.
\end{proof}

\begin{corollary}\label{cor: separable-matrix}
    Let $(X_n)_{n \in \mathbb{N}}$ and $(Z_n)_{n \in \mathbb{N}}$ be sequences of Banach spaces, $Y = \left(\bigoplus_{n \in \mathbb{N}} X_n \right)_{c_0}$, $W = \left(\bigoplus_{n \in \mathbb{N}} Z_n \right)_{c_0}$ and $S = (S_{n,m})_{n,m \in \mathbb{N}} \in \mathscr{B}(Y, W)$. Then $S$ has separable range if and only if $S_{n,m}$ has separable range for each $n, m \in \mathbb{N}$.
\end{corollary}
\begin{proof}
    One direction is elementary, so we prove the other one. Define $S_N = \sum_{n,m \leq N} \iota_n S_{n,m} \pi_m$, then $S_N \to S$ in the strong operator topology. An application of Lemma \ref{lmm: separable-to-separable} gives the result.
\end{proof}

We will also need the following elementary algebraic result.

\begin{lemma}\label{lmm: multiplicative-from-kernel}
    Let $\mathscr{B}$ and $\mathscr{C}$ be algebras over $\mathbb{K}$, and let $\mathscr{I}$ be an ideal of $\mathscr{B}$. Suppose that $\theta:\mathscr{B}\to\mathscr{C}$ is a linear map with $\ker\theta = \mathscr{I}$. If $\theta$ admits a multiplicative right inverse $\rho: \mathscr{C}\to\mathscr{B}$, then $\theta$ is multiplicative.
\end{lemma}
\begin{proof}
    Since $\rho$ is a right inverse to $\theta$, for any $b\in\mathscr{B}$ we have $b - \rho(\theta(b)) \in \ker \theta = \mathscr{I}$, so that every $b\in\mathscr{B}$ can be decomposed as $b = \rho(\theta(b)) + u$ for some $u \in \mathscr{I}$.
    
    Let $b_1,b_2 \in \mathscr{B}$ and write $b_i = \rho(\theta(b_i)) + u_i$ for some $u_i \in \mathscr{I}$, $i = 1,2$. Expanding the product,
    \begin{equation*}
    b_1 b_2 = \rho(\theta(b_1)) \rho(\theta(b_2))
    + u_1 \rho(\theta(b_2))
    + \rho(\theta(b_1)) u_2
    + u_1 u_2.
    \end{equation*}
    Since $\mathscr{I} = \ker \theta$ is an ideal the last three terms belong to $\ker\theta$, hence the linearity of $\theta$ implies
    \begin{equation*}
    \theta(b_1 b_2) = \theta(\rho(\theta(b_1)) \rho(\theta(b_2))).
    \end{equation*}
    Finally, since $\rho$ is multiplicative and a right inverse of $\theta$,
    \begin{equation*}
    \theta(\rho(\theta(b_1)) \rho(\theta(b_2))) = \theta(\rho(\theta(b_1)\theta(b_2))) = \theta(b_1)\theta(b_2),
    \end{equation*}
    which shows that $\theta$ is multiplicative.
\end{proof}

With the preceding preparations complete, we are ready to prove our next lemma.

\begin{lemma}\label{lmm: scalar-and-separable-decomposition}
    Let $(X_t)_{t \in \mathcal{P}^*}$ be a family of Banach spaces with few operators relative to $\mathcal{P}^*$, and set $Y = \left(\bigoplus_{t \in \mathcal{P}} X_{t} \right)_{c_0}$. Then for any operator $T = (T_{s,t})_{s,t \in \mathcal{P}} \in \mathscr{B}(Y)$, the matrix $\Lambda = (\lambda_{s,t})_{s,t \in \mathcal{P}}$, where $T_{s,t} = \lambda_{s,t} \iota_{t,s} + S_{s,t}$ with $\lambda_{s,t} \in \mathbb{K}$ and $S_{s,t}$ has separable range, corresponds to an operator $\theta(T) \in \incidence(\mathcal{P}^{\opposite}, c_0)$. 
    
    Moreover, the map
    \begin{equation*}
        \theta: \mathscr{B}(Y) \to \incidence(\mathcal{P}^{\opposite}, c_0), \hspace{5pt} T \mapsto \theta(T)
    \end{equation*}
    is an algebra homomorphism with $\norm{\theta} = 1$ and $\ker \theta = \mathscr{X}(Y)$.
\end{lemma}
\begin{proof}
         We start by showing that $\Lambda$ is the scalar matrix of an operator $\theta(T) \in \mathscr{B}(c_0(\mathcal{P}))$. Once this is done, it is immediate that $\theta(T) \in \incidence(\mathcal{P}^{\opposite}, c_0)$, since the family $(X_t)_{t \in \mathcal{P}}$ has few operators relative to $\mathcal{P}$.

    According to Proposition \ref{prop: when-matrix-is-operator} \ref{it: when-matrix-is-operator-b}, we need to show that $\Lambda$ has $\ell_1$-bounded rows and $c_0$-columns. Let
    \begin{equation*}
        Z = \overline{\spn}\{\iota_{s,\infty} \circ S_{s,t}[X_t]: s,t \in \mathcal{P}\},
    \end{equation*}
    which is a separable closed subspace of $X_\infty$. Since $X_0$ is non-separable and $\iota_{0,\infty}$ is an isometry, we may choose $y_0 \in B_{X_0}$ such that
    \begin{equation*}
        \delta = \dist(\iota_{0,\infty}y_0,Z)>0.
    \end{equation*}

    We first show that $\Lambda$ has $c_0$-columns. Fix $t \in \mathcal{P}$ and set $\hat{y}=\jmath_t\iota_{0,t}y_0 \in B_Y$. Compatibility of the embeddings and the definition of $Z$ give
    \begin{equation*}
        \delta|\lambda_{s,t}| \leq \norm{\iota_{s,\infty}T_{s,t}\iota_{0,t}y_0} = \norm{\pi_sT\hat{y}} \qquad(s \in \mathcal{P}).
    \end{equation*}
    Since $T\hat{y} \in Y$, the right-hand side tends to zero as $s \to \infty$. Thus $\Lambda$ has $c_0$-columns.

    We now show that $\Lambda$ has $\ell_1$-bounded rows. Fix $s \in \mathcal{P}$ and a finite set $\mathcal{T} \subseteq \mathcal{P}$. For each $t \in \mathcal{T}$, choose a unimodular scalar $\beta_{s,t}$ such that $\beta_{s,t}\lambda_{s,t}=|\lambda_{s,t}|$, and set
    \begin{equation*}
        a=\sum_{t \in \mathcal{T}}|\lambda_{s,t}|, \qquad R=\sum_{t \in \mathcal{T}}\beta_{s,t}S_{s,t}\circ\iota_{0,t}.
    \end{equation*}
    Then $R \colon X_0 \to X_s$ has separable range. Moreover, for every $z \in X_0$,
    \begin{equation*}
        \norm{a\iota_{0,s}z+Rz} = \norm{\pi_sT\left(\sum_{t \in \mathcal{T}}\beta_{s,t}\jmath_t\iota_{0,t}z\right)} \leq \norm{T}\norm{z}.
    \end{equation*}
    If $a>\norm{T}$, the reverse triangle inequality would give
    \begin{equation*}
        \norm{Rz}\geq(a-\norm{T})\norm{z} \qquad(z \in X_0).
    \end{equation*}
    Thus $R$ would be bounded below, contradicting the non-separability of $X_0$ and the separability of its range. Consequently $a\leq\norm{T}$. Since $s$ and $\mathcal{T}$ were arbitrary, we obtain
    \begin{equation*}
        \sup_{s \in \mathcal{P}}\sum_{t \in \mathcal{P}}|\lambda_{s,t}| \leq \norm{T}.
    \end{equation*}
    Hence $\Lambda$ has $\ell_1$-bounded rows and defines $\theta(T) \in \incidence(\mathcal{P}^{\opposite},c_0)$, as desired.

    Linearity of $\theta$ is automatic, and the preceding estimate gives
    \begin{equation*}
        \norm{\theta(T)} = \sup_{s \in \mathcal{P}}\sum_{t \in \mathcal{P}}|\lambda_{s,t}| \leq \norm{T}.
    \end{equation*}
    Thus $\norm{\theta}\leq1$. Since $\rho$ is an isometry and $\theta\circ\rho$ is the identity map on $\incidence(\mathcal{P}^{\opposite},c_0)$, we also have $\norm{\theta}\geq1$, so $\norm{\theta}=1$.

    By Corollary \ref{cor: separable-matrix}, we have $\ker\theta=\mathscr{X}(Y)$, and $\rho$ is multiplicative by Lemma \ref{lmm: from-scalar-going-up}. Therefore, Lemma \ref{lmm: multiplicative-from-kernel} implies that $\theta$ is multiplicative and thus a homomorphism. This finishes the proof.
\end{proof}

We can now give an analogue to Proposition \ref{prop: create-incidence-algebras} for the countable case; the proof is automatic, combining Lemmas \ref{lmm: from-scalar-going-up} and \ref{lmm: scalar-and-separable-decomposition}.

\begin{proposition}\label{prop: countable-incidence-algebras}
    Let $\mathcal{P}$ be a countable partially ordered set and let $(X_t)_{t \in \mathcal{P}^*}$ be a family of Banach spaces with few operators relative to $\mathcal{P}^*$. Set
    \begin{equation*}
        Y = \left(\bigoplus_{t \in \mathcal{P}} X_t\right)_{c_0}.
    \end{equation*}
    Then there is a strongly split short exact sequence
    \begin{equation*}
        \{0\} \longrightarrow \mathscr{X}(Y) \longrightarrow \mathscr{B}(Y) \xlongrightarrow{\theta} \incidence(\mathcal{P}^{\opposite}, c_0) \longrightarrow \{0\}.
    \end{equation*}
\end{proposition}

\begin{rem}
    The previous proposition, together with Theorem~\ref{th: almost-disjoint-families-few-operators}, enables the construction of many subalgebras of the infinite-dimensional triangular matrices as quotients of $\mathscr{B}(X)$, thereby providing a partial answer to \cite[Open Problems~(i)]{loy1989continuity}. It also yields algebras of upper and lower-triangular matrices as quotients of $\mathscr{B}(X)$ for suitable Banach spaces $X$, offering alternative solutions to \cite[Open Problems~(iv)]{loy1989continuity}, distinct from the one obtained in \cite{dales1994homomorphisms}.
\end{rem}

\section{Proof of Theorem \ref{th: quotient-algebras}.}\label{sec: proof-theorem-2}

We are finally ready to give the proof of Theorem \ref{th: quotient-algebras}.

\begin{proof}[Proof of Theorem \ref{th: quotient-algebras}]
    Let $(\mathscr{C}_n)_{n \in \mathbb{N}} \subseteq \mathfrak{A}$ be a sequence of algebras. 
    If $\mathscr C_n=\{0\}$ for every $n$, take $\mathcal A=\emptyset$. Then $C_0(K_\mathcal A)=c_0$ is separable, so $\mathscr X(c_0)=\mathscr B(c_0)$, and the zero homomorphism gives the required strongly split exact sequence. We may therefore suppose that at least one $\mathscr C_n$ is non-zero.
    For each $n \in \mathbb{N}$, we construct a partially ordered set $\mathcal{P}_n$, and throughout we assume that the sets $(\mathcal{P}_n)_{n=1}^{\infty}$ are constructed to be pairwise disjoint.
    \begin{enumerate}[label = (\roman*), ref = (\roman*)]
        \item If $\mathscr C_n=\{0\}$, we let $\mathcal P_n=\emptyset$.
        \item If $\mathscr{C}_n = \mathbb{M}_{N_n}$ for some $N_n \in \mathbb{N}$,
              or $\mathscr{C}_n = \mathscr{B}(c_0)$,
              we let $\mathcal{P}_n$ be the trivial poset (a singleton).
        \item If $\mathscr{C}_n = \incidence(\mathcal{Q}_n)$ for some finite partially ordered set $\mathcal{Q}_n$, or $\mathscr{C}_n = \incidence(\mathcal{Q}_n, c_0)$ for some countable partially ordered set $\mathcal{Q}_n$, we let $\mathcal{P}_n = \mathcal{Q}_n^{\opposite}$.
    \end{enumerate}
    
    Let $\mathcal Q=\bigcup_{n\in\mathbb N}\mathcal P_n$, ordered by
    \begin{equation*}
        s\leq_{\mathcal Q}t\quad\Longleftrightarrow\quad
        \text{there is }n\in\mathbb N\text{ such that }s,t\in\mathcal P_n
        \text{ and }s\leq_{\mathcal P_n}t.
    \end{equation*}
    Thus elements belonging to different $\mathcal P_n$ are incomparable. Define $\mathcal P=\mathcal Q^*$ to be the bounded extension from Definition~\ref{def: bounded-extension}; in particular, the two adjoined elements are respectively below and above every member of $\mathcal Q$.

    By Theorem \ref{th: almost-disjoint-families-few-operators}, we can find almost disjoint families $(\mathcal{A}_{p})_{p \in \mathcal{P}}$ such that \break $(C_0(K_{\mathcal{A}_{p}}))_{p \in \mathcal{P}}$ has few operators relative to $\mathcal{P}$.
    
    In the case where $\mathcal{P}_n$ is a singleton, we denote by $\mathcal{A}_n$ the almost disjoint family associated to it in the previous construction. If $\mathcal{P}_n$ is a poset with more than one element, its associated almost disjoint families are $(\mathcal{A}_{t})_{t \in \mathcal{P}_n}$.

    Using this, for each $n \in \mathbb{N}$, we will construct a Banach space $Y_n$, so that $Y = \left(\bigoplus_{n \in \mathbb{N}} Y_n \right)_{c_0}$ will have the desired form; again, we do this construction depending on the form of $\mathscr{C}_n$.

    \begin{enumerate}[label = (\roman*), ref = (\roman*)]
        \item \label{it: case0}If $\mathscr C_n=\{0\}$, we let $Y_n=\{0\}$.
        \item \label{it: case1} If $\mathscr{C}_n = \mathbb{M}_{N_n}$ for some $N_n \in \mathbb{N}$, we let $Y_n = C_0(K_{\mathcal{A}_n})^{N_n}$. 
        \item \label{it: case2} If $\mathscr{C}_n = \incidence(\mathcal{Q}_n)$ for some finite partially ordered set $\mathcal{Q}_n$ we \break let $Y_n = \bigoplus_{t \in \mathcal{P}_n} C_0(K_{\mathcal{A}_{t}})$.
        \item \label{it: case4} If $\mathscr{C}_n = \incidence(\mathcal{Q}_n, c_0)$ for some countable partially ordered set $\mathcal{Q}_n$ we let $Y_n = \left(\bigoplus_{t \in \mathcal{P}_n} C_0(K_{\mathcal{A}_{t}}) \right)_{c_0}$.
        \item \label{it: case3} If $\mathscr{C}_n = \mathscr{B}(c_0)$, we let
        \begin{equation*}
            Y_n = \left(C_0(K_{\mathcal{A}_n}) \oplus C_0(K_{\mathcal{A}_n}) \oplus \ldots \right)_{c_0} = c_0(C_0(K_{\mathcal{A}_n})) = C_0(\omega \times K_{\mathcal{A}_n}).
        \end{equation*}
    \end{enumerate}

    It follows that for every $n \in \mathbb{N}$, we have a strongly split exact sequence
    \begin{equation}\label{eq: split-exact-sequence}
        \{0\} \longrightarrow \mathscr{X}\left(Y_n \right) \longrightarrow \mathscr{B}\left(Y_n \right)
            \xlongrightarrow{\theta_n} \mathscr{C}_n \longrightarrow \{0\}
    \end{equation}
    with an isometric right inverse $\rho_n$ and $\norm{\theta_n}\leq1$, with equality whenever $\mathscr C_n\neq\{0\}$. For case~\ref{it: case0} this is the zero sequence. The assertion is immediate for case~\ref{it: case1}; for case~\ref{it: case2}, it follows from Proposition~\ref{prop: create-incidence-algebras}; for case~\ref{it: case4}, it follows from Proposition~\ref{prop: countable-incidence-algebras}; and for case~\ref{it: case3}, it follows from \cite[Remark 1.7]{acuaviva2025operators} (alternatively, by arguing as in the $c_0$-incidence algebra case).

    Let $Y = \left(\bigoplus_{n \in \mathbb{N}} Y_n \right)_{c_0}$, so that by Lemma \ref{lmm: combine-pairwise-disjoint-families}, $Y \cong C_0(K_\mathcal{A})$ for some almost disjoint family $\mathcal{A} \subseteq \inset{\mathbb{N}}$. We claim that $\mathcal{A}$ satisfies Theorem \ref{th: quotient-algebras}.  \\
    
    Recall that any operator $T: Y \to Y$ can be naturally identified with an infinite matrix of operators $(T_{n,m})_{n,m\in \mathbb{N}}$ where $T_{n,m} \in \mathscr{B}(Y_m, Y_n)$. Elementary arguments show that the map $\rho: \prod_{n \in \mathbb{N}} \mathscr{C}_n \to \mathscr{B}(Y)$, defined by mapping $a = (a_n)_{n \in \mathbb{N}} \in \prod_{n \in \mathbb{N}} \mathscr{C}_n$ to the operator $\rho(a)$ given by the infinite diagonal matrix $(\delta_{n,m} \rho_n(a_n))_{n,m \in \mathbb{N}}$ is a well-defined isometric homomorphism.

    Let $T: Y \to Y$ with $T = (T_{n,m})_{n,m \in \mathbb{N}}$. By \eqref{eq: split-exact-sequence}, it follows that for each $n \in \mathbb{N}$, $T_{n,n} = \rho_n(\theta_n(T_{n,n})) + S_{n,n}$ where $S_{n,n}$ has separable range. Moreover, from the construction and Corollary \ref{cor: separable-matrix}, $T_{n,m}$ has separable range whenever $n \not = m$. Observe that the sequence $(\theta_n(T_{n,n}))_{n \in \mathbb{N}}$ is uniformly bounded since both $\theta_n$ and $T_{n,n}$ are; thus $(\theta_n(T_{n,n}))_{n \in \mathbb{N}} \in \prod_{n \in \mathbb{N}} \mathscr{C}_n$. 
    
    It is clear that the map $\theta: \mathscr{B}(Y) \to \prod_{n \in \mathbb{N}} \mathscr{C}_n$, $T \mapsto (\theta_n(T_{n,n}))_{n \in \mathbb{N}}$ is a bounded linear map with $\norm{\theta} \leq 1$. Since $\rho$ is a right inverse to $\theta$, it follows that $\norm{\theta} = 1$. Moreover, from the previous discussion and Corollary \ref{cor: separable-matrix}, we have $\ker \theta = \mathscr{X}(Y)$. By Lemma \ref{lmm: multiplicative-from-kernel}, it follows that $\theta$ is multiplicative and thus an algebra homomorphism. This gives the desired strongly split exact sequence, which completes the proof.
\end{proof}

\section{Applications}\label{sec:applications}

This section has two parts. We first use the quotient construction to obtain a discontinuous homomorphism from an operator algebra onto the upper-triangular matrix algebra $\mathbb T_2(\mathbb K)$. We then examine the range of possible quotient algebras and use their ideal lattices to count the closed ideals of the corresponding operator algebras.

\subsection{A discontinuous homomorphism}\label{subsec:zfc-discontinuous}

The proof proceeds in two stages. First, we construct a discontinuous surjective homomorphism from a suitable $c_0$-incidence algebra onto $\mathbb T_2(\mathbb K)$. We then apply Theorem~\ref{th: quotient-algebras} to realise this incidence algebra as a strongly split quotient of an algebra $\mathscr B(C_0(K_{\mathcal A}))$. The bounded right inverse supplied by that theorem ensures that discontinuity is preserved when the homomorphism is composed with the quotient map. A final extension argument passes from $C_0(K_{\mathcal A})$ to $C(K)$.

\subsubsection{Construction of the discontinuous homomorphism} The main difficulty in the first stage is to construct an unbounded functional on the off-diagonal block which is nevertheless compatible with multiplication by the two diagonal blocks. For this purpose, we construct finite-to-one maps $f_n\colon D\to R$ and ultrafilters $\mathcal U$ and $\mathcal V$ such that all the pushforwards $f_n(\mathcal U)$ equal $\mathcal V$, while distinct maps $f_n$ and $f_m$ disagree on a set belonging to $\mathcal U$. The resulting ultrafilter matrix coefficients $F_n$ separate a sequence of norm-one operators $S_n$. Passing to an algebraic quotient and using a Hamel-basis extension then produces an unbounded functional $L$ satisfying the required bimodule identities. These identities are precisely what makes the resulting map into $\mathbb T_2(\mathbb K)$ multiplicative.

We write $\mathbb T_2(\mathbb K)$ for the algebra of matrices
\begin{equation*}
    \begin{pmatrix}
        \alpha & \beta\\
        0 & \gamma
    \end{pmatrix}
    \qquad (\alpha,\beta,\gamma\in\mathbb K),
\end{equation*}
equipped with the operator norm inherited from $\mathscr B(\ell_\infty^2)$.

If $\mathcal U$ is an ultrafilter on a set $S$, define $\chi_{\mathcal U}\colon\ell_\infty(S)\to\mathbb K$ by
\begin{equation*}
    \chi_{\mathcal U}(a)=\lim_{s\to\mathcal U}a(s)
    \qquad (a\in\ell_\infty(S)).
\end{equation*}
The limit is well defined because a closed disk in $\mathbb K$ is compact, and continuity of the scalar operations shows that $\chi_{\mathcal U}$ is a norm-one character. Namely, $\chi_{\mathcal U}$ denotes the character given by ultrafilter evaluation. If $\mathcal U$ is principal at $s_0$, then $\chi_{\mathcal U}(a)=a(s_0)$, while a free ultrafilter gives a generalised limiting value. For a map $f\colon S\to T$, we write
\begin{equation*}
    f(\mathcal U)=\{A\subseteq T:f^{-1}(A)\in\mathcal U\}
\end{equation*}
for the pushforward ultrafilter.

The defining property of the pushforward gives the change-of-variables identity
\begin{equation}\label{eq:pushforward-ultralimit}
    \lim_{s\to\mathcal U}b(f(s))=\lim_{t\to f(\mathcal U)}b(t)\qquad(b\in\ell_\infty(T)),
\end{equation}
which will be the reason for synchronising the pushforwards of the maps used below.

For each integer $k\geq1$, choose a finite set $I_k$ so that the sets $I_k$ are pairwise disjoint and $|I_k|\geq k2^k$. Set
\begin{equation*}
    R=\bigsqcup_{k\geq1}I_k,
    \qquad
    D_k=\operatorname{Inj}(\{1,\ldots,k\},I_k),
    \qquad
    D=\bigsqcup_{k\geq1}D_k,
\end{equation*}
where $\operatorname{Inj}(A,B)$ denotes the set of injective maps from $A$ to $B$.

An element $x\in D_k$ may be regarded as an ordered list of $k$ distinct elements of $I_k$. Define $f_n\colon D\to R$ by
\begin{equation*}
f_n(x)=
\begin{cases}
x(n),&x\in D_k\text{ and }n\leq k,\\
x(1),&x\in D_k\text{ and }n>k.
\end{cases}
\end{equation*}
That is, the map $f_n$ reads the $n$th entry of $x$ whenever that entry exists, and otherwise reads its first entry. For each integer $n\geq1$, each $f_n$ is finite-to-one: $f_n(D_j)\subseteq I_j$ for every $j\geq1$, and the sets $I_j$ are pairwise disjoint. Hence, if $r\in I_k$, then $f_n^{-1}(\{r\})\subseteq D_k$, which is finite.

The next lemma shows that the maps $f_n$ can have a common pushforward along an ultrafilter $\mathcal U$, while remaining pairwise distinct $\mathcal U$-almost everywhere.

\begin{lemma}\label{lmm:synchronised-ultrafilters}
There are free ultrafilters $\mathcal U$ on $D$ and $\mathcal V$ on $R$ such that
\begin{equation*}
    f_n(\mathcal U)=\mathcal V
    \qquad (n\geq1)
\end{equation*}
and
\begin{equation*}
    \{x\in D:f_n(x)=f_m(x)\}\notin\mathcal U
    \qquad (n,m\geq1,\ n\neq m).
\end{equation*}
\end{lemma}

\begin{proof}
We construct $\mathcal U$ by specifying a family of subsets of $D$ that it must avoid.

For integers $n,m\geq1$ and $A\subseteq R$, set
\begin{equation*}
    E(n,m,A)
    =f_n^{-1}(A)\mathbin{\triangle}f_m^{-1}(A),
\end{equation*}
where $\triangle$ denotes symmetric difference. For $K\geq1$, set
\begin{equation*}
    B_K=D_1\sqcup\ldots\sqcup D_K.
\end{equation*}
The sets $E(n,m,A)$ encode the requirement that the pushforwards of the maps $f_n$ coincide, while the sets $B_K$ encode the requirement that $\mathcal U$ be free.

Let $\mathcal J$ consist of all subsets of finite unions of sets of the form $E(n,m,A)$ and $B_K$. Equivalently, $\mathcal J$ is the ideal of the Boolean algebra $\wp(D)$ generated by these sets. Here an ideal of subsets means a family closed under taking subsets and finite unions. We claim that $\mathcal J$ is proper, that is, $D\notin\mathcal J$.
We show that no finite union of generators of $\mathcal J$ covers $D$. Since the sets $(B_K)_{K\geq1}$ are increasing, it is enough to fix some $K, r\geq1$ and consider sets of the form
\begin{equation*}
    B_K\cup\bigcup_{j=1}^r E(n_j,m_j,A_j).
\end{equation*}
Choose an integer
\begin{equation*}
    k>\max\{K,r,n_1,m_1,\ldots,n_r,m_r\}.
\end{equation*}
For each $\varepsilon=(\varepsilon_1,\ldots,\varepsilon_r)\in\{0,1\}^r$, define
\begin{equation*}
    C_\varepsilon=\{t\in I_k:1_{A_j}(t)=\varepsilon_j\text{ for every }1\leq j\leq r\},
\end{equation*}
where $1_{A_j}$ denotes the indicator function of $A_j$. The sets $C_\varepsilon$ partition $I_k$ into at most $2^r$ classes. Hence, by the pigeonhole principle, there is some $\varepsilon\in\{0,1\}^r$ such that
\begin{equation*}
    |C_\varepsilon|\geq\frac{|I_k|}{2^r}\geq k2^{k-r}\geq k.
\end{equation*}

Choose distinct elements $c_1,\ldots,c_k\in C_\varepsilon$, and define $x\colon\{1,\ldots,k\}\to I_k$ by $x(i)=c_i$. Then $x\in D_k$. Moreover, since $k>K$, we have $x\notin B_K$.

Fix $1\leq j\leq r$. Since $n_j,m_j<k$, we have $f_{n_j}(x)=c_{n_j}$ and $f_{m_j}(x)=c_{m_j}$. Both points belong to $C_\varepsilon$, so
\begin{equation*}
    1_{A_j}(c_{n_j})=\varepsilon_j=1_{A_j}(c_{m_j}).
\end{equation*}
Consequently,
\begin{equation*}
    f_{n_j}(x)\in A_j\quad\Longleftrightarrow\quad f_{m_j}(x)\in A_j,
\end{equation*}
and therefore $x\notin E(n_j,m_j,A_j)$.

Thus
\begin{equation*}
    x\notin B_K\cup\bigcup_{j=1}^r E(n_j,m_j,A_j).
\end{equation*}
Since the finite collection of generators was arbitrary, no finite union of generators covers $D$. Hence $D\notin\mathcal J$, and $\mathcal J$ is proper.

Since $\mathcal J$ is proper, the family
\begin{equation*}
    \mathcal F=\{D\setminus J:J\in\mathcal J\}
\end{equation*}
is a proper filter on $D$. Indeed, closure of $\mathcal J$ under finite unions makes $\mathcal F$ closed under finite intersections, closure of $\mathcal J$ under taking subsets makes $\mathcal F$ closed under taking supersets, and properness of $\mathcal J$ ensures that $\emptyset\notin\mathcal F$. By the ultrafilter lemma, extend $\mathcal F$ to an ultrafilter $\mathcal U$ on $D$. Since $B_K\in\mathcal J$, we have $D\setminus B_K\in\mathcal U$ for every $K\geq1$. Moreover, every finite subset of $D$ is contained in some $B_K$. Therefore no finite subset of $D$ belongs to $\mathcal U$, and hence $\mathcal U$ is free.

For every $A\subseteq R$ and $n,m\geq1$, the complement of $E(n,m,A)$ belongs to $\mathcal U$. It follows that
\begin{equation*}
    f_n^{-1}(A)\in\mathcal U
    \quad\Longleftrightarrow\quad
    f_m^{-1}(A)\in\mathcal U.
\end{equation*}
Indeed, two sets whose symmetric difference has complement in an ultrafilter must have the same membership in that ultrafilter. Thus all the pushforwards $f_n(\mathcal U)$ coincide; denote their common value by $\mathcal V$. To see that $\mathcal V$ is free, fix $n\geq1$ and suppose that a finite set $F\subseteq R$ belonged to $\mathcal V=f_n(\mathcal U)$. Then $f_n^{-1}(F)\in\mathcal U$ by the definition of pushforward. Since $f_n$ is finite-to-one, this preimage is finite, contradicting the freeness of $\mathcal U$.

Finally, suppose that $n\neq m$. If $x\in D_k$ with $k\geq\max\{n,m\}$, then the injectivity of $x$ gives $f_n(x)\neq f_m(x)$. Hence
\begin{equation*}
    \{x\in D:f_n(x)=f_m(x)\}
    \subseteq B_{\max\{n,m\}-1}.
\end{equation*}
The agreement set is therefore finite, so it does not belong to the free ultrafilter $\mathcal U$.
\end{proof}

Set
\begin{equation*}
    \mathscr M=\mathscr B(c_0(R),c_0(D)).
\end{equation*}
For $Q\in\{D,R\}$ and $h\in\ell_\infty(Q)$, let $M_h^Q\in\mathscr B(c_0(Q))$ denote the multiplication operator defined by
\begin{equation*}
    (M_h^Qz)(q)=h(q)z(q)\qquad(q\in Q).
\end{equation*}
For each integer $n\geq1$, define
\begin{equation*}
    F_n(T)=\lim_{x\to\mathcal U}e_x^*Te_{f_n(x)}\qquad(T\in\mathscr M).
\end{equation*}
The scalar function $x\mapsto e_x^*Te_{f_n(x)}$ is bounded in modulus by $\norm{T}$, so $F_n$ is a well-defined linear functional with $\norm{F_n}\leq1$. Thus $F_n$ reads the matrix entry of $T$ in row $x$ and column $f_n(x)$ and then takes its $\mathcal U$-limit. Moreover, for $a\in\ell_\infty(D)$ and $b\in\ell_\infty(R)$, we have
\begin{equation}\label{eq:bimodule-Fn}
    F_n(M_a^D T M_b^R)=\chi_{\mathcal U}(a)F_n(T)\chi_{\mathcal V}(b)
\end{equation}

Indeed, the two multiplication operators act on the relevant coordinate vectors by
\begin{equation*}
    M_b^R e_{f_n(x)}=b(f_n(x))e_{f_n(x)},\qquad e_x^*M_a^D=a(x)e_x^*.
\end{equation*}
Consequently,
\begin{equation*}
    e_x^*M_a^D T M_b^R e_{f_n(x)}=a(x)\bigl(e_x^*T e_{f_n(x)}\bigr)b(f_n(x)),
\end{equation*}
Taking the $\mathcal U$-limit separates the three scalar factors. The first two limits are $\chi_{\mathcal U}(a)$ and $F_n(T)$, while~\eqref{eq:pushforward-ultralimit} and $f_n(\mathcal U)=\mathcal V$ give
\begin{equation*}
    \lim_{x\to\mathcal U}b(f_n(x))=\lim_{r\to\mathcal V}b(r)=\chi_{\mathcal V}(b).
\end{equation*}
This proves~\eqref{eq:bimodule-Fn}.

For each integer $n\geq1$, define $S_n\colon c_0(R)\to c_0(D)$ by
\begin{equation*}
    (S_nz)(x)=z(f_n(x))\qquad(x\in D).
\end{equation*}

If $z\in c_0(R)$ and $\varepsilon>0$, then
\begin{equation*}
    \{x\in D:|(S_nz)(x)|\geq\varepsilon\}=f_n^{-1}(\{r\in R:|z(r)|\geq\varepsilon\}),
\end{equation*}
The set on the right is finite: the level set $\{r\in R:|z(r)|\geq\varepsilon\}$ is finite because $z\in c_0(R)$, and its preimage is finite because $f_n$ is finite-to-one. Hence $S_nz\in c_0(D)$ and $S_n\in\mathscr M$. It is immediate that $\norm{S_n}\leq1$, and equality follows by taking $x_0\in D$ and $z=e_{f_n(x_0)}$. The choice of $S_n$ allows the functionals $F_m$ to distinguish the operators $S_n$. Indeed, since $S_nz=z \circ f_n$, for each $x\in D$ we have
\begin{equation*}
    e_x^*S_ne_{f_m(x)} =(S_ne_{f_m(x)})(x) =e_{f_m(x)}(f_n(x)) =1_{\{y\in D:f_n(y)=f_m(y)\}}(x).
\end{equation*}
It follows from the definition of $F_m$ that
\begin{equation*}
    F_m(S_n)
    =\lim_{x\to\mathcal U}
      1_{\{y\in D:f_n(y)=f_m(y)\}}(x).
\end{equation*}
If $m=n$, the function inside the limit is identically one. If $m\neq n$, Lemma~\ref{lmm:synchronised-ultrafilters} shows that its support does not belong to $\mathcal U$, and hence its $\mathcal U$-limit is zero. Therefore
\begin{equation}\label{eq:Fn-Sn}
    F_m(S_n)=\delta_{m,n}\qquad(m,n\geq1).
\end{equation}

We now convert this bounded biorthogonal system into an unbounded functional that satisfies bimodule relations with respect to the characters $\chi_\mathcal{U}$ and $\chi_\mathcal{V}$.

\begin{lemma}\label{lmm:unbounded-bimodule-functional}
There is an unbounded linear functional $L\colon\mathscr M\to\mathbb K$ such that
\begin{equation*}
    L(M_a^D T)=\chi_{\mathcal U}(a)L(T),\qquad L(T M_b^R)=L(T)\chi_{\mathcal V}(b)
\end{equation*}
for all $T\in\mathscr M$, $a\in\ell_\infty(D)$, and $b\in\ell_\infty(R)$.
\end{lemma}

\begin{proof}
Let $\mathscr N$ be the algebraic subspace
\begin{equation*}
    \mathscr N=\spn\{M_a^D T:\chi_{\mathcal U}(a)=0,\ T\in\mathscr M\}+\spn\{T M_b^R:\chi_{\mathcal V}(b)=0,\ T\in\mathscr M\}.
\end{equation*}
No closure is taken. Taking $b=1_R$ in~\eqref{eq:bimodule-Fn} gives
\begin{equation*}
    F_n(M_a^D T)=\chi_{\mathcal U}(a)F_n(T),
\end{equation*}
so this expression vanishes whenever $\chi_{\mathcal U}(a)=0$. Similarly, taking $a=1_D$ gives
\begin{equation*}
    F_n(T M_b^R)=F_n(T)\chi_{\mathcal V}(b),
\end{equation*}
which vanishes whenever $\chi_{\mathcal V}(b)=0$. By linearity, every $F_n$ vanishes on $\mathscr N$.

We show that the cosets
\begin{equation*}
    S_n+\mathscr N\qquad(n\geq1)
\end{equation*}
are linearly independent in $\mathscr M/\mathscr N$.

Indeed, suppose that for some $r\geq1$ and $c_1,\ldots,c_r\in\mathbb K$,
\begin{equation*}
    \sum_{j=1}^r c_j(S_j+\mathscr N)=0
\end{equation*}
in $\mathscr M/\mathscr N$. Then $\sum_{j=1}^r c_jS_j\in\mathscr N$, and hence, for each $1\leq m\leq r$,
\begin{equation*}
    0=F_m\left(\sum_{j=1}^r c_jS_j\right)=\sum_{j=1}^r c_jF_m(S_j)=c_m
\end{equation*}
by~\eqref{eq:Fn-Sn}; this shows the desired linear independence.

We may therefore define a linear functional on their linear span by
\begin{equation*}
    \lambda_0\left(\sum_{j=1}^r c_j(S_j+\mathscr N)\right)=\sum_{j=1}^r jc_j.
\end{equation*}
Extend the family $\{S_n+\mathscr N:n\geq1\}$ to a Hamel basis of $\mathscr M/\mathscr N$ over $\mathbb K$. Assigning the value zero to the additional basis vectors extends $\lambda_0$ to a linear functional
\begin{equation*}
    \lambda\colon\mathscr M/\mathscr N\to\mathbb K.
\end{equation*}
Let $q\colon\mathscr M\to\mathscr M/\mathscr N$ be the canonical algebraic quotient map and put $L=\lambda\circ q$. Then
\begin{equation*}
    L(S_n)=n,\qquad\norm{S_n}=1\qquad(n\geq1),
\end{equation*}
so $L$ is unbounded.

For $a\in\ell_\infty(D)$ and $T\in\mathscr M$, we have $\chi_{\mathcal U}(a-\chi_{\mathcal U}(a)1_D)=0$, and therefore
\begin{equation*}
    M_a^D T-\chi_{\mathcal U}(a)T=M_{a-\chi_{\mathcal U}(a)1_D}^D T\in\mathscr N.
\end{equation*}
Similarly, for $b\in\ell_\infty(R)$,
\begin{equation*}
    T M_b^R-\chi_{\mathcal V}(b)T=T M_{b-\chi_{\mathcal V}(b)1_R}^R\in\mathscr N.
\end{equation*}
Since $\mathscr N=\ker q\subseteq\ker L$, we get
\begin{equation*}
    L(M_a^D T)=\chi_{\mathcal U}(a)L(T),\qquad L(T M_b^R)=L(T)\chi_{\mathcal V}(b),
\end{equation*}
which shows the desired bimodule relation and finishes the proof.
\end{proof}

\subsubsection{The incidence algebra} It remains to assemble the two diagonal characters and the off-diagonal functional $L$ into a single homomorphism. For this, we use the height-two poset in which every element of $D$ lies below every element of $R$. Its $c_0$-incidence algebra has exactly the required $2\times2$ block form: multiplication operators on the diagonal and the full operator space $\mathscr M$ in the upper-right corner.

Let $\mathcal P=D\sqcup R$, ordered by retaining equality and declaring
\begin{equation*}
    d<r\qquad(d\in D,\ r\in R).
\end{equation*}
This is a countable poset. Relative to $c_0(\mathcal P)=c_0(D)\oplus_\infty c_0(R)$, Definition~\ref{def: E-inc-algebra} gives
\begin{equation}\label{eq:bipartite-block-form}
    \incidence(\mathcal P,c_0)=\left\{\begin{pmatrix}M_a^D&T\\0&M_b^R\end{pmatrix}:a\in\ell_\infty(D),\ b\in\ell_\infty(R),\ T\in\mathscr M\right\}.
\end{equation}
Indeed, distinct elements of $D$ are incomparable, as are distinct elements of $R$, so both diagonal blocks are multiplication operators. The block from $c_0(D)$ to $c_0(R)$ must vanish, while the block from $c_0(R)$ to $c_0(D)$ is unrestricted.

Define
\begin{equation*}
    \Phi\colon\incidence(\mathcal P,c_0)\to\mathbb T_2(\mathbb K),\qquad \Phi\left(\begin{pmatrix}M_a^D&T\\0&M_b^R\end{pmatrix}\right)=\begin{pmatrix}\chi_{\mathcal U}(a)&L(T)\\0&\chi_{\mathcal V}(b)\end{pmatrix}.
\end{equation*}
This is an algebra homomorphism. Indeed, let
\begin{equation*}
    A=\begin{pmatrix}M_a^D&T\\0&M_b^R\end{pmatrix},\qquad B=\begin{pmatrix}M_{a'}^D&S\\0&M_{b'}^R\end{pmatrix}.
\end{equation*}
The diagonal entries of $\Phi(AB)$ agree with those of $\Phi(A)\Phi(B)$ because $\chi_{\mathcal U}$ and $\chi_{\mathcal V}$ are characters. The upper-right block of $AB$ is $M_a^D S+T M_{b'}^R$, and Lemma~\ref{lmm:unbounded-bimodule-functional} gives
\begin{equation*}
    L(M_a^D S+T M_{b'}^R)=\chi_{\mathcal U}(a)L(S)+L(T)\chi_{\mathcal V}(b'),
\end{equation*}
which is the upper-right entry of $\Phi(A)\Phi(B)$.

The map $\Phi$ is surjective because, for all $\alpha,\beta,\gamma\in\mathbb K$,
\begin{equation*}
    \Phi\left(\begin{pmatrix}M_{\alpha 1_D}^D&\beta S_1\\0&M_{\gamma 1_R}^R\end{pmatrix}\right)=\begin{pmatrix}\alpha&\beta\\0&\gamma\end{pmatrix}.
\end{equation*}
It is discontinuous because
\begin{equation*}
    A_n=\begin{pmatrix}0&S_n\\0&0\end{pmatrix}\qquad(n\geq1)
\end{equation*}
satisfies $\norm{A_n}=1$, whereas
\begin{equation*}
    \Phi(A_n)=\begin{pmatrix}0&n\\0&0\end{pmatrix}
\end{equation*}
has norm $n$.

We shall also need the following elementary observation.

\begin{lemma}\label{lmm:one-point-compactification}
Let $K_0$ be a locally compact Hausdorff space, let $K$ be its one-point compactification, and let $\mathscr C$ be a Banach algebra. Suppose that
\begin{equation*}
    \Psi\colon\mathscr B(C_0(K_0))\to\mathscr C
\end{equation*}
is an algebra homomorphism whose kernel contains every finite-rank operator on $C_0(K_0)$. Then there is an algebra homomorphism
\begin{equation*}
    \widetilde\Psi\colon\mathscr B(C(K))\to\mathscr C
\end{equation*}
having the same range as $\Psi$. Moreover, $\widetilde\Psi$ is continuous if and only if $\Psi$ is continuous, and it is compact if and only if $\Psi$ is compact. In particular, surjectivity and discontinuity are preserved by this lifting.
\end{lemma}

\begin{proof}
Identify $C_0(K_0)$ with the closed subspace of functions in $C(K)$ that vanish at $\infty$. Let $\jmath\colon C_0(K_0)\to C(K)$ be the canonical inclusion, let $u=1_K$, and define $\varepsilon\in C(K)^*$ and $\pi\colon C(K)\to C_0(K_0)$ by
\begin{equation*}
    \varepsilon(f)=f(\infty),\qquad\pi(f)=f-\varepsilon(f)u\qquad(f\in C(K)).
\end{equation*}
Using the vector-first convention for rank-one operators, we have $\pi\jmath=I_{C_0(K_0)}$ and
\begin{equation*}
    I_{C(K)}=\jmath\pi+u\otimes\varepsilon.
\end{equation*}
Define
\begin{equation*}
    \widetilde\Psi(T)=\Psi(\pi T\jmath)\qquad(T\in\mathscr B(C(K))).
\end{equation*}
For $T,S\in\mathscr B(C(K))$, we have
\begin{equation*}
    \pi TS\jmath=(\pi T\jmath)(\pi S\jmath)+(\pi Tu)\otimes(\varepsilon S\jmath).
\end{equation*}
The second summand has rank at most one and therefore belongs to $\ker\Psi$. Hence $\widetilde\Psi(TS)=\widetilde\Psi(T)\widetilde\Psi(S)$.

Define
\begin{equation*}
    \iota\colon\mathscr B(C_0(K_0))\to\mathscr B(C(K)),\qquad\iota(A)=\jmath A\pi.
\end{equation*}
Since $\pi\jmath=I_{C_0(K_0)}$, the map $\iota$ is a bounded algebra homomorphism and
\begin{equation*}
    \widetilde\Psi(\iota(A))=\Psi(A)\qquad(A\in\mathscr B(C_0(K_0))).
\end{equation*}
Since every value of $\widetilde\Psi$ is also a value of $\Psi$, the two homomorphisms have the same range.

The map $T\mapsto\pi T\jmath$ is bounded, so continuity and compactness of $\Psi$ imply the corresponding properties of $\widetilde\Psi$. Conversely, $\Psi=\widetilde\Psi\circ\iota$, and $\iota$ is bounded, so continuity and compactness of $\widetilde\Psi$ imply the corresponding properties of $\Psi$.
\end{proof}

We are now ready for the construction of the desired discontinuous homomorphism. 

\begin{proof}[Proof of Theorem~\ref{th:zfc-discontinuous-homomorphism}]
Apply Theorem~\ref{th: quotient-algebras} to the sequence
\begin{equation*}
    \mathscr C_1=\incidence(\mathcal P,c_0),\qquad\mathscr C_n=\{0\}\qquad(n\in\mathbb N\setminus\{1\}).
\end{equation*}
We obtain an almost disjoint family $\mathcal A$, a bounded surjective algebra homomorphism
\begin{equation*}
    \theta\colon\mathscr B(C_0(K_{\mathcal A}))\longrightarrow\incidence(\mathcal P,c_0),
\end{equation*}
and an isometric algebra-homomorphic right inverse
\begin{equation*}
    \rho\colon\incidence(\mathcal P,c_0)\longrightarrow\mathscr B(C_0(K_{\mathcal A}))
\end{equation*}
such that $\theta\circ\rho$ is the identity on $\incidence(\mathcal P,c_0)$. Moreover,
\begin{equation*}
    \ker\theta=\mathscr X(C_0(K_{\mathcal A})).
\end{equation*}

Set $\Psi=\Phi\circ\theta$. Since both $\Phi$ and $\theta$ are surjective, so is $\Psi$. If $\Psi$ were continuous, then $ \Phi=\Psi\circ\rho$ would be continuous, contradicting the construction of $\Phi$. Thus $\Psi$ is a discontinuous surjective algebra homomorphism from $\mathscr B(C_0(K_{\mathcal A}))$ onto $\mathbb T_2(\mathbb K)$.

Let $K$ be the one-point compactification of $K_{\mathcal A}$. Since
\begin{equation*}
    \mathscr X(C_0(K_{\mathcal A}))=\ker\theta\subseteq\ker\Psi,
\end{equation*}
every finite-rank operator on $C_0(K_{\mathcal A})$ belongs to $\ker\Psi$. Lemma~\ref{lmm:one-point-compactification} therefore gives an algebra homomorphism
\begin{equation*}
    \widetilde\Psi\colon\mathscr B(C(K))\longrightarrow\mathbb T_2(\mathbb K)
\end{equation*}
having the same range as $\Psi$. Hence $\widetilde\Psi$ is surjective, and it is discontinuous because continuity is preserved in both directions by the lifting.

Finally, the points $x_n$ indexed by $\mathbb N$ are isolated in $K$, and
\begin{equation*}
    K'\subseteq\{y_A:A\in\mathcal A\}\cup\{\infty\}.
\end{equation*}
Each $y_A$ is isolated in $K'$, so $K''\subseteq\{\infty\}$ and $K'''=\emptyset$. Thus $K$ is compact, scattered, and has Cantor--Bendixson height at most three.
\end{proof}

\subsection{Possible quotient algebras and ideal counts}\label{subsec:possible-quotient-algebras}

Theorem~\ref{th: quotient-algebras} gives a large supply of quotient algebras of the form $\mathscr B(C_0(K_{\mathcal A}))/\mathscr X(C_0(K_{\mathcal A}))$. In this final subsection, we determine the finite ideal counts obtained from the building blocks in that theorem and give some concrete examples. We also show that the construction can produce the largest possible number of maximal closed ideals.

If $\mathscr C$ is a Banach algebra, let $\mathfrak I(\mathscr C)$ denote the lattice of closed two-sided ideals of $\mathscr C$, including $\{0\}$ and $\mathscr C$ itself. For a finite partially ordered set $\mathcal P$, write
\begin{equation*}
    \operatorname{Int}(\mathcal P)=\{(s,t)\in\mathcal P\times\mathcal P:s\leq t\},
\end{equation*}
ordered by containment of intervals,
\begin{equation*}
    (s,t)\preceq(u,v)\quad\Longleftrightarrow\quad u\leq s\leq t\leq v.
\end{equation*}
Set
\begin{equation*}
    N(\mathcal P)=|\operatorname{Up}(\operatorname{Int}(\mathcal P))|,
\end{equation*}
where $\operatorname{Up}(\mathcal Q)$ denotes the family of upsets of a partially ordered set $\mathcal Q$. By \cite[Theorem~3.1]{doubilet1972foundations}, the lattice of two-sided ideals of $\incidence(\mathcal P)$ is naturally anti-isomorphic to the lattice of order ideals of $\operatorname{Int}(\mathcal P)$. Since $\mathcal P$ is finite, all these ideals are closed, and since complementation identifies order ideals with upsets,
\begin{equation*}
    N(\mathcal P)=|\mathfrak I(\incidence(\mathcal P))|.
\end{equation*}

\begin{proposition}\label{prop:finite-ideal-counts}
The finite numbers of non-trivial closed ideals furnished by Theorem~\ref{th: quotient-algebras} are precisely the integers of the form
\begin{equation*}
    3^r\prod_{i=1}^kN(\mathcal P_i),
\end{equation*}
where $\mathcal P_1,\ldots,\mathcal P_k$ are finite partially ordered sets and $r,k \geq 0$. If $r+k\geq1$, then, for each such choice, the corresponding algebra $\mathscr B(C_0(K_{\mathcal A}))$ may be chosen to have exactly
\begin{equation*}
    r+\sum_{i=1}^k|\mathcal P_i|
\end{equation*}
maximal closed ideals. In the exceptional case $r=k=0$, the corresponding algebra may be chosen to be $\mathscr B(c_0)$, which has one non-trivial closed ideal, and this ideal is also its unique maximal closed ideal.
\end{proposition}

\begin{proof}
If $r=k=0$, apply Theorem~\ref{th: quotient-algebras} to the identically zero sequence. As in the proof of that theorem, we may take $\mathcal A=\emptyset$, so that $C_0(K_{\mathcal A})=c_0$. The algebra $\mathscr B(c_0)$ has precisely three closed ideals,
\begin{equation*}
    \{0\},\qquad\mathscr K(c_0),\qquad\mathscr B(c_0).
\end{equation*}

Thus it has exactly one non-trivial closed ideal, which is also its unique maximal closed ideal. We may therefore suppose that $r+k\geq1$.

First, apply Theorem~\ref{th: quotient-algebras} to a sequence whose non-zero terms are
\begin{equation*}
    \incidence(\mathcal P_1),\ldots,\incidence(\mathcal P_k)
\end{equation*}
together with $r$ copies of $\mathscr B(c_0)$. We obtain an almost disjoint family $\mathcal A$ such that
\begin{equation*}
    \mathscr B(C_0(K_{\mathcal A}))/\mathscr X(C_0(K_{\mathcal A}))\cong\left(\prod_{i=1}^k\incidence(\mathcal P_i)\right)\times\mathscr B(c_0)^r.
\end{equation*}
The closed ideals of a finite product of unital Banach algebras are precisely the products of closed ideals in the coordinates. Each algebra $\incidence(\mathcal P_i)$ has $N(\mathcal P_i)$ closed ideals, while $\mathscr B(c_0)$ has exactly three closed ideals \cite{gohberg1967normally}. Hence the quotient has
\begin{equation*}
    3^r\prod_{i=1}^kN(\mathcal P_i)
\end{equation*}
closed ideals. Proposition~\ref{prop: quotient-to-ideals} gives the same number of non-trivial closed ideals in $\mathscr B(C_0(K_{\mathcal A}))$. Indeed, passing from the ideals of the quotient to the non-trivial ideals of $\mathscr B(C_0(K_{\mathcal A}))$ removes the whole algebra and adds the missing ideal $\mathscr K(C_0(K_{\mathcal A}))$, so the number is unchanged.

Conversely, suppose that
\begin{equation*}
    \prod_{n\in\mathbb N}\mathscr C_n
\end{equation*}
is a quotient algebra furnished by Theorem~\ref{th: quotient-algebras} and has only finitely many closed ideals. If this product is zero, then its lattice of closed ideals has one element, and Proposition~\ref{prop: quotient-to-ideals} gives exactly one non-trivial closed ideal in the corresponding operator algebra. We may therefore suppose that the quotient is non-zero.

Since the quotient has only finitely many closed ideals, its non-zero factors are matrix algebras, finite incidence algebras, and copies of $\mathscr B(c_0)$; an infinite $c_0$-incidence algebra is excluded because the diagonal characters $T\mapsto e_t^*Te_t$ have pairwise distinct maximal kernels. Each matrix algebra has exactly two closed ideals, a contribution already represented by a one-point poset. Since the number of closed ideals of a finite product is the product of the numbers of closed ideals of its factors, the stated form follows.

For the maximal-ideal count, \cite[Corollary~3.2]{doubilet1972foundations} shows that for a finite poset $\mathcal{P}$ the maximal ideals of $\incidence(\mathcal P)$ are precisely
\begin{equation*}
    \{f\in\incidence(\mathcal P):f(t,t)=0\}\qquad(t\in\mathcal P).
\end{equation*}
Thus $\incidence(\mathcal P)$ has exactly $|\mathcal P|$ maximal ideals. The algebra $\mathscr B(c_0)$ has a unique maximal closed ideal, and a maximal ideal of a finite product is obtained by choosing a maximal ideal in one coordinate and taking the whole algebra in every other coordinate. Hence the quotient constructed above has
\begin{equation*}
    r+\sum_{i=1}^k|\mathcal P_i|
\end{equation*}
maximal ideals. Proposition~\ref{prop: quotient-to-ideals} transfers these maximal ideals to $\mathscr B(C_0(K_{\mathcal A}))$.
\end{proof}

It remains to identify concrete values of $N(\mathcal P)$; the product formula then supplies arbitrary finite products of these values, together with arbitrary powers of $3$. Let us spell out some examples.

If $\mathcal P$ is an antichain with $d$ points, then $\operatorname{Int}(\mathcal P)$ is an antichain with $d$ points, so
\begin{equation*}
    N(\mathcal P)=2^d.
\end{equation*}
Thus every power of $2$ is obtained. If $\mathcal P$ is a chain with $d$ points, then the upsets of $\operatorname{Int}(\mathcal P)$ are counted by the Catalan number $C_{d+1}$, giving
\begin{equation*}
    2,5,14,42,132,\ldots.
\end{equation*}

The values arising from height-two posets admit a useful general formula. Let $G=(L\sqcup R,E)$ be a finite bipartite graph, where $L\sqcup R$ stands for the vertex set of the graph and $E$ for its edges. Let $\mathcal P_G$ be the poset on the disjoint union $L\sqcup R$ defined by
\begin{equation*}
    \ell<r\quad\Longleftrightarrow\quad\ell\in L,\ r\in R,\ \text{and }(\ell,r)\in E.
\end{equation*}
For $I\subseteq L$ and $J\subseteq R$, write
\begin{equation*}
    e_G(I,J)=|E\cap(I\times J)|.
\end{equation*}

\begin{proposition}\label{prop: bipartite-height-two-count}
For every finite bipartite graph $G=(L\sqcup R,E)$,
\begin{equation*}
    N(\mathcal P_G)=\sum_{I\subseteq L}\sum_{J\subseteq R}2^{e_G(I,J)}.
\end{equation*}
\end{proposition}

\begin{proof}
The elements of $\operatorname{Int}(\mathcal P_G)$ are the diagonal intervals $(v,v)$, for $v\in L\sqcup R$, together with the intervals $(\ell,r)$ corresponding to the edges $(\ell,r)\in E$. The only non-trivial order relations are
\begin{equation*}
    (\ell,\ell)\preceq(\ell,r)\qquad\text{and}\qquad(r,r)\preceq(\ell,r).
\end{equation*}
Let $\mathcal W$ be an upset, and let $I\subseteq L$ and $J\subseteq R$ consist of those vertices $v$ for which $(v,v)\notin\mathcal W$. Once $I$ and $J$ are fixed, every edge interval having at least one endpoint outside $I\sqcup J$ is forced to belong to $\mathcal W$. The remaining edge intervals correspond to the $e_G(I,J)$ edges between $I$ and $J$, and each may be included independently. Thus there are $2^{e_G(I,J)}$ upsets with these prescribed missing diagonal intervals. Summing over $I$ and $J$ proves the formula.
\end{proof}

If $G=K_{p,q}$ is complete bipartite, Proposition~\ref{prop: bipartite-height-two-count} gives
\begin{equation*}
    N(\mathcal P_{K_{p,q}})=\sum_{i=0}^{p}\sum_{j=0}^{q}\binom pi\binom qj2^{ij}=\sum_{i=0}^{p}\binom pi(1+2^i)^q.
\end{equation*}
For $p=1$, this gives
\begin{equation*}
    N(\mathcal P_{K_{1,q}})=2^q+3^q
\end{equation*}
and hence we obtain the values
\begin{equation*}
    5,13,35,97,\ldots.
\end{equation*}
For $p=2$, we have
\begin{equation*}
    N(\mathcal P_{K_{2,q}})=2^q+2\cdot3^q+5^q,
\end{equation*}
which yields
\begin{equation*}
    47,187,803,\ldots
\end{equation*}
for $q=2,3,4,\ldots$. Similar constructions yield further values, and arbitrary finite products of these values are also obtainable. However, this does not appear to cover all the natural numbers, so Question~\ref{quest:numberideals} remains open.

We observe that we can also obtain the largest possible number of maximal closed ideals for an algebra of cardinality $\mathfrak c$.

\begin{corollary}\label{cor: two-to-c-maximal-ideals}
There is an almost disjoint family $\mathcal A\subseteq\inset{\mathbb N}$ such that $\mathscr B(C_0(K_{\mathcal A}))$ has exactly $2^{\mathfrak c}$ maximal closed ideals.
\end{corollary}

\begin{proof}
Apply Theorem~\ref{th: quotient-algebras} with $\mathscr C_n=\mathbb K$ for every $n$. The resulting quotient is
\begin{equation*}
    \ell_\infty=C(\beta\mathbb N).
\end{equation*}
Its maximal ideals are the kernels of evaluation at the points of $\beta\mathbb N$, and $|\beta\mathbb N|=2^{\mathfrak c}$. Proposition~\ref{prop: quotient-to-ideals} therefore gives the same number of maximal closed ideals in $\mathscr B(C_0(K_{\mathcal A}))$.
\end{proof}

Let us finally indicate the scope of the construction. Theorem~\ref{th: quotient-algebras} realises countable normed products of the listed building blocks, but it does not characterize all Banach algebras that can occur as such quotients. Accordingly, Proposition~\ref{prop:finite-ideal-counts} describes the ideal counts furnished by this construction, rather than all possible counts arising from arbitrary almost disjoint families. \\

\noindent\textbf{Acknowledgements.} This paper forms part of the author’s PhD research at Lancaster University, conducted under the supervision of Professor N. J. Laustsen. The author is deeply grateful to Professor Laustsen for his insightful comments and valuable suggestions, which significantly improved the presentation of the manuscript. He also wishes to thank both Professor Laustsen and Dr. Horv\'ath for kindly explaining the ideas behind Lemma~\ref{lmm: niels-and-bence}, and for granting permission to include it in this work. Finally, the author thanks Witold Marciszewski for pointing out that the proof of Theorem~\ref{th: almost-disjoint-families-few-operators} extends naturally from partially ordered sets of cardinality at most $\mathfrak c$ to partially ordered sets which order-embed into $\wp(\mathfrak c)$.

He also acknowledges with thanks the funding from the EPSRC (grant number EP/W524438/1) that has supported his studies. \\

For the purpose of open access, the author has applied a Creative Commons Attribution (CC BY) licence to any Author Accepted Manuscript version arising. \medskip

\noindent\textbf{Data availability.} No data was used for the research described in the article.


\end{document}